\documentclass[a4paper, 11pt]{article}

\usepackage{latexsym}
\usepackage{amscd}
\usepackage{amsmath}
\usepackage{amssymb}

\usepackage{amsthm}
\usepackage{float}
\usepackage{graphicx}
\usepackage{psfrag}
\usepackage{pst-node}
\usepackage{pst-plot}

\theoremstyle{plain}
\newtheorem{theorem}{Theorem}
\newtheorem{corollary}[theorem]{Corollary}
\newtheorem{lemma}[theorem]{Lemma}

\newtheorem{proposition}[theorem]{Proposition}

\newtheorem{step}{Step}

\theoremstyle{definition}
\newtheorem{definition}[theorem]{Definition}
\newtheorem{remark}[theorem]{Remark}

\newdimen\argwidth
\def\db[#1\db]{%
  \setbox0=\hbox{$#1$}\argwidth=\wd0
  \setbox0=\hbox{$\left[\box0\right]$}
    \advance\argwidth by -\wd0
  \left[\kern.3\argwidth\box0 \kern.3\argwidth\right]}

\newcommand{\diag}{\ensuremath{\mathrm{diag}}}

\renewcommand{\Im}{\ensuremath{\mathop{\mathfrak{Im}}\nolimits}}
\newcommand{\Ker}{\mathrm{Ker}}

\newcommand{\Coker}{\mathrm{Coker}}
\newcommand{\Ob}{\mathfrak{Ob}}

\newcommand{\id}{\ensuremath{\mathop{\mathrm{id}}}}

\newcommand{\Ext}{\mathop{\mathrm{Ext}}\nolimits}
\newcommand{\Hom}{\mathop{\mathrm{Hom}}\nolimits}
\newcommand{\RHom}{\mathop{\bR\mathrm{Hom}}\nolimits}

\newcommand{\Lotimes}{\stackrel{\bL}{\otimes}}
\newcommand{\End}{\mathop{\mathrm{End}}\nolimits}
\newcommand{\Spec}{\operatorname{Spec}}

\newcommand{\SL}{SL}
\newcommand{\bCx}{\bC^{\times}}

\newcommand{\coh}{\operatorname{coh}}

\newcommand{\module}{\operatorname{mod}}
\newcommand{\Per}[1]{\!{\ }^{#1} \! \mathop{\mathrm{Per}}}

\newcommand{\Fuk}{\mathop{\mathfrak{Fuk}}\nolimits}
\newcommand{\m}{\mathfrak{m}}

\newcommand{\bA}{\ensuremath{\mathbb{A}}}

\newcommand{\bC}{\ensuremath{\mathbb{C}}}

\newcommand{\bF}{\ensuremath{\mathbb{F}}}

\newcommand{\bL}{\ensuremath{\mathbb{L}}}

\newcommand{\bQ}{\ensuremath{\mathbb{Q}}}
\newcommand{\bR}{\ensuremath{\mathbb{R}}}

\newcommand{\bZ}{\ensuremath{\mathbb{Z}}}

\newcommand{\scA}{\ensuremath{\mathcal{A}}}
\newcommand{\scB}{\ensuremath{\mathcal{B}}}
\newcommand{\scC}{\ensuremath{\mathcal{C}}}
\newcommand{\scD}{\ensuremath{\mathcal{D}}}
\newcommand{\scE}{\ensuremath{\mathcal{E}}}
\newcommand{\scF}{\ensuremath{\mathcal{F}}}

\newcommand{\scH}{\ensuremath{\mathcal{H}}}

\newcommand{\scL}{\ensuremath{\mathcal{L}}}
\newcommand{\scM}{\ensuremath{\mathcal{M}}}

\newcommand{\scO}{\ensuremath{\mathcal{O}}}
\newcommand{\scP}{\ensuremath{\mathcal{P}}}

\newcommand{\scT}{\ensuremath{\mathcal{T}}}

\newcommand{\scW}{\ensuremath{\mathcal{W}}}

\newcommand{\scZ}{\ensuremath{\mathcal{Z}}}

\newcommand{\frakS}{\ensuremath{\mathfrak{S}}}

\newcommand{\frakh}{\ensuremath{\mathfrak{h}}}

\newcommand{\reg}{\mathrm{reg}}
\newcommand{\BG}{B_n} % braid group
\newcommand{\ABG}{B_n^{(1)}} % affine braid group
\newcommand{\Diff}{\operatorname{Diff}}

\newcommand{\palg}{\scA_n^{(1)}}

\newcommand{\Stab}{\operatorname{Stab}}
\newcommand{\Supp}{\ensuremath{\operatorname{Supp}}}

  \newcommand{\Auteq}{\ensuremath{\operatorname{Auteq}}}
  
  \newcommand{\Span}[1]{\left<#1\right>}
\newcommand{\Br}{\operatorname{Br}}
\newcommand{\Aut}{\operatorname{Aut}}

\newcommand{\GLt}{\widetilde{GL^+}}
\newcommand{\wt}{\widetilde}

\title{Stability Conditions on $A_n$-Singularities}
\author{Akira Ishii, Kazushi Ueda, and Hokuto Uehara}
\date{}
\begin{document}

\maketitle

\begin{abstract}
We study the spaces of locally finite
stability conditions
on the derived categories of coherent sheaves
on the minimal resolutions of $A_n$-singularities
supported at the exceptional sets.
Our main theorem is that
they are connected
and
simply-connected.
The proof
is based on the study of
spherical objects
in \cite{Ishii-Uehara_ADC}
and
the homological mirror symmetry
for $A_n$-singularities.
\end{abstract}

%%%%%%%%%%%%%%%%%%%%%%%%%%%%%%%%%%%%%%%%%%%%%%%%%%%%%%%%%%%%%%%%%%%%%%%%
%%%%%%%%%%%%%%%%%%%%%%%%%%%%%%%%%%%%%%%%%%%%%%%%%%%%%%%%%%%%%%%%%%%%%%%%
\section{Introduction}
The theory of
stability conditions on triangulated categories
is introduced by Bridgeland
\cite{Bridgeland_SCTC}
based on the work of
Douglas et al.
\cite{Aspinwall-Douglas_DSM,
Douglas_DCNS,
Douglas_DCYM,
Douglas_DBHMSS,
Douglas-Fiol-Roemelsberger_SBB}
on the stability of BPS D-branes.
It is a fine mixture
of the theory of $t$-structures
\cite{Beilinson-Bernstein-Deligne}
and the slope stability
\cite{Mumford_ICM62},
which allows us to represent
any object in a triangulated category
as a successive mapping cone of
semistable objects
in a unique way.
He proved that
the set
$\Stab \scT$
of stability conditions
on a triangulated category $\scT$
satisfying an additional assumption
called local-finiteness
has a natural structure of a complex manifold,
and proposed to study this manifold
as an invariant of $\scT$.
Since the definition of stability conditions
uses only the triangulated structure of $\scT$,
the group $\Auteq \scT$
of triangle autoequivalences of $\scT$
naturally acts on the manifold $\Stab \scT$,
suggesting a geometric approach
to study
%geometric properties of $\Stab \scT$ reflects
the structure of $\Auteq \scT$.

%This problem of using $\Stab \scT$
%to get a grip of $\Auteq \scT$
Such an approach
has been pursued by Bridgeland himself
\cite{Bridgeland_SCKS}
when
$\scT$ is
the bounded derived category
$D^b \coh X$
of coherent sheaves
on a complex algebraic $K3$ surface $X$,
leading him to the following
remarkable result and conjecture:
There is a distinguished connected component
$\Sigma(X)$
of the space $\Stab X$
of locally finite, numerical stability conditions on $D^b \coh X$.
His conjecture is:
\begin{enumerate}
 \item $\Sigma(X)$ is preserved by $\Auteq D^b \coh X$, and
 \item $\Sigma(X)$ is simply-connected.
\end{enumerate}
Assuming this conjecture,
he could prove that
$\Auteq D^b \coh X$ is an extension of
the index two subgroup
$\Aut^{+} H^*(X, \bZ)$
of the group of Hodge isometries of the Mukai lattice
of $X$
by the fundamental group
$\pi_1 \scP_0^+(X)$
of the period domain of $X$.

For
%a field $k$ and
a positive integer $n$,
let
$$
 f : X \to Y = \Spec \bC[x, y, z] / (x y + z^{n+1})
$$
be the minimal resolution of the $A_n$-singularity.
Let further $\scD$ be
the bounded derived category $D^b \coh_Z X$
of coherent sheaves on $X$
supported at the exceptional set $Z$,
and $\scC$ be
its full triangulated subcategory
consisting of objects $E$
satisfying $\bR f_* E = 0$.
The categories
$\scC$ and $\scD$ serve
as toy models of
the derived categories of coherent sheaves
on $K3$ surfaces.
The main result in this paper
is the following:

\begin{theorem} \label{theorem:main}
$\Stab \scC$ and $\Stab \scD$ are connected,
and
$\Stab \scD$ is simply-connected.
\end{theorem}

This result,
together with
the simply-connectedness of
a distinguished connected component of $\Stab \scC$
proved by Thomas \cite{Thomas_SCBG},
shows that
the above conjecture of Bridgeland
holds in these cases.
When $n=1$,
the connectedness of $\Stab \scD$
has also been proved by
Okada \cite{Okada_SMCYS}.

The basic strategy of our proof
for the connectedness of $\Stab \scD$
is to find
a stability condition
such that
structure sheaves of all the closed points
are stable
in any given connected component of $\Stab \scD$.
Since the set of such stability conditions
form a distinguished connected open subset of $\Stab \scD$,
the connectedness of $\Stab \scD$ follows.

Theorem \ref{theorem:covering} due to Bridgeland
shows that
the simply-connectedness of $\Stab \scD$
follows from the faithfulness of
an affine braid group action on $\scD$.
We prove this
using homological mirror symmetry
for $A_n$-singularities
and ideas from
Khovanov and Seidel
\cite{Khovanov-Seidel}.
Unfortunately,
we have to work over a field of characteristic two
in order to apply Theorem \ref{theorem:Khovanov-Seidel}
by Khovanov and Seidel,
and we lift this faithfulness result
to any characteristic
using the deformation theory of complexes
by Inaba \cite{Inaba_TDMCCS}.

In contrast to the case of $\Stab \scD$,
we cannot use algebro-geometric argument
in the proof of the connectedness of $\Stab \scC$,
since $\scC$ does not contain
any skyscraper sheaves.
Instead, we use
a result in \cite{Ishii-Uehara_ADC}
and ideas from \cite{Khovanov-Seidel}
to reduce the problem
of the connectedness of $\Stab \scC$
to that of configurations of curves
on a disk.

The organization of this paper is as follows:
In \S \ref{section:generalities},
we collect basic definitions and known results
used in this paper.
In \S \ref{section:McKay},
we recall the McKay correspondence
for $A_n$-singularities
in such a way
that is valid in any characteristic.
In \S \ref{section:connectedness_D},
we give the proof
of the connectedness of $\Stab \scD$.
We prove the faithfulness in characteristic two
in \S \ref{section:char_two}
and lift it to any characteristic
in \S \ref{section:lifting}.
The connectedness of $\Stab \scC$ is proved
in \S \ref{section:connectedness_C}.
In the appendix,
we prove that every autoequivalence of $\scD$
is given by an integral functor.

{\bf Acknowledgment}:
We thank Jun-ichi Matsuzawa and
Yukinobu Toda for valuable discussions and suggestions.
A.~I.~ is supported by the Grants-in-Aid 
for Scientific Research (No.18540034).
K.~U.~is supported
by the 21st Century COE Program
of Osaka University.
H.~U.~is supported by the Grants-in-Aid 
for Scientific Research (No.17740012).
 
%%%%%%%%%%%%%%%%%%%%%%%%%%%%%%%%%%%%%%%%%%%%%%%%%%%%%%%%%%%%%%%%%%%%%%%%
%%%%%%%%%%%%%%%%%%%%%%%%%%%%%%%%%%%%%%%%%%%%%%%%%%%%%%%%%%%%%%%%%%%%%%%%
\section{Generalities} \label{section:generalities}

We collect basic definitions
and known results
in this section.
All the categories
appearing in this paper
will be essentially small.
For a triangulated category $\scT$,
$K(\scT)$ denotes its Grothendieck group,
and for an object $E \in \scT$,
$[E]$ will denote its class
in $K(\scT)$.
For two objects $E, F \in \scT$
and $i \in \bZ$,
$\Hom^*_{\scT}(E, F)$,
$\Hom^{\leq i}_{\scT}(E, F)$, and
$\Hom^{\geq i}_{\scT}(E, F)$ will denote
$\bigoplus_{j \in \bZ} \Hom^j_{\scT}(E, F)$,
$\bigoplus_{j \leq i} \Hom^j_{\scT}(E, F)$, and
$\bigoplus_{j \geq i} \Hom^j_{\scT}(E, F)$, respectively.

%%%%%%%%%%%%%%%%%%%%%%%%%%%%%%%%%%%%%%%%%%%%%%%%%%%%%%%%%%%
\subsection{Stability conditions on triangulated categories}

The following definition is introduced
by Bridgeland \cite{Bridgeland_SCTC}
based on the work of
Douglas et al.
\cite{Aspinwall-Douglas_DSM,
Douglas_DCNS,
Douglas_DCYM,
Douglas_DBHMSS,
Douglas-Fiol-Roemelsberger_SBB}
on the stability of BPS D-branes:

\begin{definition}
A stability condition
$\sigma = (Z,\scP)$
on a triangulated category $\scT$
consists of
\begin{itemize}
 \item a group homomorphism $Z:K(\scT)\to\bC$, and
 \item full additive subcategories
       $\scP(\phi)$ for $\phi\in\bR$
\end{itemize}
satisfying the following conditions:
\begin{enumerate}
\item
If $0\ne E\in \scP(\phi)$, then $Z(E)=m(E)\exp(i\pi\phi)$ for some
$m(E)\in\bR_{>0}$,
\item
for all $\phi\in\bR$, $\scP(\phi+1)=\scP(\phi)[1]$,
\item
for $A_j\in\scP(\phi_j)$ $(j=1,2)$ with $\phi_1>\phi_2$, we have
$\Hom_{\scT}(A_1,A_2)=0$,
\item
for every non-zero object $E\in\scT$, there is a finite sequence of
real numbers
\[\phi_1>\phi_2>\cdots>\phi_n\]
and a collection of triangles
\begin{equation} \label{eq:HN_filtration}
\psmatrix[colsep=.5, rowsep=.7]
  0 = \Rnode{a}{E_0} & & E_1 & & E_2 & \cdots & E_{n-1} & &  
\Rnode{n}{E_n}=E \\
        & A_1 & & A_2 &   &        &       & A_n & \\
\endpsmatrix
\psset{nodesep=3pt,arrows=->}
\ncline{a}{1,3}
\ncline{1,3}{1,5}
\ncline{1,6}{1,7}
\ncline{1,7}{n}
\ncline{1,3}{2,2}
\ncline{1,5}{2,4}
\ncline{n}{2,8}
\psset{arrows=-}
\ncline{1,5}{1,6}
\psset{arrows=->,linestyle=dashed}
\ncline{2,2}{a}
\ncline{2,4}{1,3}
\ncline{2,8}{1,7}
\end{equation}
with $A_j\in\scP(\phi_j)$ for all $j$.
\end{enumerate}
\end{definition}

$Z$ is called the {\em central charge},
and the collection of triangles
in (\ref{eq:HN_filtration})
is called the {\em Harder-Narasimhan filtration}.
It follows from the definition
that $\scP(\phi)$ is an abelian category,
and its non-zero object $\scE \in \scP(\phi)$ is said to be
{\em semistable of phase $\phi$}.
$\scE$ is said to be {\em stable}
if it is a simple object of $\scP(\phi)$,
i.e., there are no proper subobjects
of $\scE$
in $P(\phi)$.
By \cite[proposition 5.3]{Bridgeland_SCTC},
to give a stability condition
on a triangulated category $\scT$
is equivalent to
giving a bounded $t$-structure on $\scT$
and a {\em stability function}
(previously called a centered slope-function)
on its heart
with the {\em Harder-Narasimhan property}.
For the definitions of
a stability function
and the Harder-Narasimhan property,
see \cite[\S 2]{Bridgeland_SCTC}.

The set of stability conditions
satisfying a certain technical condition
called \emph{local-finiteness}
\cite[definition 5.7]{Bridgeland_SCTC}
is denoted by $\Stab \scT$.
This condition ensures
that each $\scP(\phi)$ is a finite length category
so that each semi-stable
object has a Jordan-H{\" o}lder filtration.
By combining it with the Harder-Narasimhan filtration,
any non-zero object $E \in \scT$
admits a decomposition as
in (\ref{eq:HN_filtration})
such that
$A_j \in \scP(\phi_j)$ is stable
for all $j$
and
$\phi_1 \ge \phi_2 \ge \dots \ge \phi_n$.
Bridgeland introduces a natural topology
on $\Stab\scT$
such that the forgetful map
$$
\begin{array}{rccc}
 \scZ : & \Stab \scT & \to & \Hom(K(\scT), \bC) \\
 & \rotatebox{90}{$\in$} & & \rotatebox{90}{$\in$} \\
 & (Z, \scP) & \mapsto & Z
\end{array}
$$
satisfies the following:
\begin{theorem}[{\cite[theorem 1.2]{Bridgeland_SCTC}}]
\label{theorem:annals}
For each connected component $\Sigma$ of $\Stab \scT$,
there is a linear subspace
$
 V(\Sigma) \subset \Hom(K(\scT), \bC)
$
with a well-defined linear topology
such that the restriction $\scZ |_{\Sigma}$
gives a local homeomorphism.
\end{theorem}
Hence $\Stab\scT$ forms
a (possibly infinite-dimensional)
complex manifold
modeled on the topological vector space
$V(\Sigma)$.
When $\scT = \scC$ or $\scD$,
$K(\scT)$ is finite-dimensional,
and we prove in Lemma \ref{lemma:perturb}
that $V(\Sigma)$ always coincides with
$\Hom(K(\scT), \bC)$.

Since the definition of $\Stab \scT$
uses only the triangulated structure of $\scT$,
the group $\Auteq \scT$
of triangle autoequivalences of $\scT$
acts naturally on $\Stab \scT$
from the left;
for $\sigma = (Z, \scP) \in \Stab \scT$
and
$\Phi \in \Auteq \scT$,
$$
 \Phi(\sigma) = (\Phi^{*} Z, \Phi(\scP))
$$
where
$\Phi^*$ is the pull-back
by the inverse of the automorphism
$
 \Phi_* : K(\scT) \to K(\scT)
$
induced by $\Phi$.
This action commutes with the right action
of the universal cover
$\GLt(2, \bR)$
of the general linear group
$GL^+(2, \bR)$
with positive determinant,
which ``rotates'' the central charge
\cite[lemma 8.2]{Bridgeland_SCTC}.

%%%%%%%%%%%%%%%%%%%%%%%%%%%%%%%%%%%%%%%%%%%%%%%%%%%%%%%%%%%
\subsection{Minimal resolutions of $A_n$-singularities}
We consider an arbitrary field $k$.
The case $\operatorname{char}(k)=2$ will be important later.
For % a field $k$ and
a positive integer $n$,
let
$$
 f : X \to \Spec k[x, y, z] / (x y + z^{n+1})
$$
be the minimal resolution of
the $A_n$-singularity.
The exceptional set of $f$
will be denoted by
$$
 Z = f^{-1}(0) = C_1 \cup \dots \cup C_n,
$$
where $C_i$'s are irreducible $(-2)$-curves
such that
$
 C_i \cap C_j = \emptyset
$
if $|i - j| > 1$.
Let $\scD_k$ be
the bounded derived category of coherent sheaves
on $X$ supported at $Z$
and
$\scC_k$ be its full triangulated subcategory
consisting of objects $E$
satisfying $\bR f_* E = 0$.
Put
$E_0 = \omega_Z$ and
$E_i = \scO_{C_i}(-1)$ for $i = 1, \dots, n$.
Here, $\omega_Z$ is the dualizing sheaf of $Z$.
Then we have
$$
 \scC_k = \langle E_1, \dots, E_n \rangle
$$
and
$$
 \scD_k = \langle E_0, \dots, E_n \rangle,
$$
where $\langle \bullet \rangle$
denotes the smallest full triangulated subcategory
of $\scD_k$
containing them.
We simply write $\scC$ and $\scD$ instead of $\scC_{\bC}$
and $\scD_{\bC}$, respectively.

For $E, F\in \scD_k$,
define the Euler form by
\begin{equation}\label{equatyion:euler1}
\chi(E, F) = \sum_i (-1)^i \dim \Hom_{\scD_k}^i(E, F),
\end{equation}
which descends to a bilinear form
on $K(\scD_k)$.
By the Riemann-Roch formula,
we have
\[
\chi(E, F) = - c_1(E) \cdot c_1(F).
\]
The Euler form $\chi$ endows $K(\scD_k)$
with the structure of
the affine root lattice
of type $A_n^{(1)}$.
A non-zero element $\alpha \in K(\scD)$ is a {\it root} if $\chi(\alpha, \alpha) \le 2$
and it is a {\it real root} if $\chi(\alpha, \alpha)=2$.
An {\it imaginary root} is a root that is not a real root.
Let $\delta \in K(\scD_k)$ be the class of
the structure sheaf of a closed point
with residue field $k$.
\label{def:delta}
Then an imaginary root
is a non-zero element of
$\bZ \delta \subset K(\scD_k)$.

%%%
%%% {lemma:root}

\begin{lemma}\label{lemma:root}
If $E \in \scD$ is stable
with respect to
some stability condition,
then
$[E] \in K(\scD)$
is a root.
\end{lemma}

\begin{proof}
The stability of $E$ implies
$\Hom_{\scD}^{\leq -1}(E, E)=0$ and
$\Hom_{\scD}(E, E) \cong \bC$.
The Serre duality shows
$\Hom^{\geq 3}_{\scD}(E, E)=0$ and
$\Hom^2_{\scD}(E, E) \cong \bC$.
Hence $\chi(E, E) \le 2$
and $[E]$ is a root.
\end{proof}

\begin{definition}
\begin{enumerate}
 \item
An object $E \in \scD_k$
is spherical
if 
$$
 \Hom_{\scD_k}^{i}(E, E)
  \cong \begin{cases}
     k & \text{if} \ i = 0, 2, \\
     0 & \text{otherwise}.
    \end{cases}
$$
 \item
An ordered set $(E_1, \dots, E_n)$
of spherical objects in $\scD_k$
is an $A_n$-configuration
if
$$
 \dim \Hom_{\scD_k}^*(E_i, E_j)
  = \begin{cases}
     1 & \text{if} \ |i-j|=1, \\
     0 & \text{if} \ |i-j| \geq 2.
    \end{cases}
$$
\end{enumerate}
\end{definition}

The proof of Lemma \ref{lemma:root}
shows the following:

\begin{lemma} \label{lemma:stable_real_root_spherical}
If the class $[E] \in K(\scD)$
of a stable object $E \in \scD$
is a real root,
then $E$ is spherical.
\end{lemma}

A spherical object $E \in \scD_k$
gives rise to an autoequivalence of $\scD_k$
through the {\em twist functor} $T_E$,
defined as the Fourier-Mukai transform
with
$$
 \{ E^{\vee} \boxtimes E \to \scO_\Delta \}
  \in D^b \coh X \times X
%T_E(F) = \{\RHom(E, F) \otimes E \to F \}
$$
as the kernel \cite{Seidel-Thomas}.
Define
$$
 \Br(\scD_k) = \langle T_{E_0}, \dots, T_{E_n} \rangle
  \subset \Auteq \scD_k
$$
and
$$
 \Br(\scC_k) = \langle T_{E_1}, \dots, T_{E_n} \rangle
  \subset \Auteq \scC_k.
$$

Define the braid group $\BG$
as the group
generated by
$\sigma_1, \dots, \sigma_n$
subject to relations
\begin{eqnarray*}
 \sigma_i \sigma_{i+1} \sigma_i
   &=& \sigma_{i+1} \sigma_i \sigma_{i+1}, \quad i=1, \dots, n-1, \\
 \sigma_i \sigma_j &=& \sigma_j \sigma_i, \qquad \qquad \quad |i-j|>2.
\end{eqnarray*}
It has the following topological description:
Let
$$
 \frakh = \{ (a_1, \dots, a_{n+1}) \in \bC^{n+1} \mid 
               a_1 + \dots + a_{n+1} = 0 \}
$$
be a Cartan subalgebra
of the complex simple Lie algebra
of type $A_n$,
and $\frakh^\reg$
be the complement of its root hyperplanes,
$$
 \frakh^\reg = \{ (a_1, \dots, a_{n+1}) \in \frakh \mid
                  a_i \ne a_j \ \text{for} \ i \neq j \}.
$$
The Weyl group $W \cong \frakS_{n+1}$ acts on $\frakh$
by permutations,
and $\frakh^\reg$ is the set of regular orbits
of $W$.
Then $\BG$ is isomorphic
to the fundamental group
of the quotient $\frakh^\reg / W$
\cite{Brieskorn_FRROEKS, Deligne_IGTG}.
It follows that
$\BG$ has another topological description:
Let
$
 \Delta = \{ 1, \zeta, \zeta^2, \dots, \zeta^{n} \}
$
be the set of $(n+1)$th roots of unity
and
$
 \Diff_0 (\bC)
$
be the group of
diffeomorphisms of $\bC$
which are the identity map
outside compact sets.
Then there is a map
$
 \Diff_0(\bC) \to \frakh^\reg / W
$
which sends $\phi \in \Diff_0(\bC)$ to
$[\{\phi(1)-c,\phi(\zeta)-c,\dots, \phi(\zeta^n)-c\}]$ with $c=\sum_{i=0}^n \phi(\zeta^n) / (n+1)$.
This map is a Serre fibration whose fiber
over $[\Delta]$ is the subgroup
$
 \Diff_0(\bC; \Delta) \subset \Diff_0(\bC)
$
which fixes $\Delta$ as a set.
From the long exact sequence of homotopy groups
associated to this fibration,
we can see that
$$
 \BG \cong \pi_0(\Diff_0(\bC; \Delta)).
$$

The assignment
$\sigma_i \mapsto T_{E_i}$
for $i=1, \dots, n$
defines a homomorphism
from $\BG$
to $\Br(\scC)$,
which is injective
by Khovanov, Seidel, and Thomas
\cite{Khovanov-Seidel, Seidel-Thomas}.
This result is the key
to the proof of the simply-connectedness
of a distinguished connected component
of $\Stab \scC$
by Thomas \cite{Thomas_SCBG}.

Now define the affine braid group $\ABG$
to be the group
generated by $\sigma_0, \dots, \sigma_n$
subject to relations
\begin{eqnarray*}
 \sigma_i \sigma_{i+1} \sigma_i
   &=& \sigma_{i+1} \sigma_i \sigma_{i+1}, \quad i=0, \dots, n, \\
 \sigma_i \sigma_j &=& \sigma_j \sigma_i, \qquad \quad \quad |i-j|>2.
\end{eqnarray*}
Here,
we put
$\sigma_{n+1} = \sigma_0$
by notation.
Let $\widehat{\frakh}^\reg$
be the complement of the affine root hyperplanes
in $\widehat{\frakh} = \frakh \oplus \bC$:
$$
 \widehat{\frakh}^\reg
  = \{ (a_1, \dots, a_{n+1}, b) \in \frakh \oplus \bC \mid
        a_i - a_j + b d \ne 0 \ \text{for} \ 
        i \ne j \ \text{and} \ d \in \bZ \}.
$$
The affine Weyl group $\widehat{W}$
acts freely on $\widehat{\frakh}^\reg$,
and the fundamental group of the orbit space
$\widehat{\frakh}^\reg / \widehat{W}$
is given by $\ABG \times \bZ$
by \cite{Dung_FGSROAWG}.

The group $\ABG$ also admits
the following topological interpretation:
Let $\Diff_0(\bCx)$ be
the group of diffeomorphisms of $\bCx$
which are the identity maps
outside compact sets
and
$\Diff_0(\bCx; \Delta)$ be its subgroup
fixing $\Delta$ as a set.
We can define a homomorphism
$$
 \ABG \to \pi_0 ( \Diff_0 (\bCx; \Delta ) )
$$
from $\ABG$
to the group of
connected components
of $\Diff_0 (\bCx; \Delta)$
by sending $\sigma_i$
to the class of a diffeomorphism of $\bCx$
which permutes two neighboring points
$\zeta^i$ and $\zeta^{i+1}$
for $i=0, \dots, n$.
This homomorphism
is known to be injective
(cf. \cite{Kent-Peifer_AGADABG}).

The assignment
$\sigma_i \mapsto T_{E_i}$
for $i=0, \dots, n$
defines a homomorphism
$$
\rho : \ABG \to \Br(\scD_k) ,
$$
which
can be extended
to a surjective homomorphism
$$
\tilde{\rho}: \ABG \times \bZ \to \Br(\scD_k) \times \bZ.
$$
Here, 
since $\Br(\scD_k)$ does not contain
any power of the shift functor,
the right-hand side is considered
as a subgroup of $\Auteq \scD_k$
so that the second factor $\bZ$
corresponds to the group
generated by the shift functor [2].
In \cite{Bridgeland_SCKS}, Bridgeland shows 
the following theorem
for any Kleinian singularities, that is, 
rational double points over $\bC$.

\begin{theorem}[{\cite[theorem 1.3]{Bridgeland_SCKS}}]
\label{theorem:covering}
There is a connected component of $\Stab \scD$
which is a covering space of
$\widehat{\frakh}^\reg / \widehat{W}$
such that the group $\Br(\scD) \times \bZ$
acts as the group of deck transformations.
\end{theorem}

Hence we have the canonical group homomorphism 
$$\pi_1(\widehat{\frakh}^\reg / \widehat{W})=\ABG \times \bZ
\to \Br(\scD) \times \bZ.$$
Bridgeland also shows that
this homomorphism coincides with $\tilde{\rho}$.
Therefore if $\rho$ is injective, we conclude 
that the connected component in 
Theorem \ref{theorem:covering} 
is simply connected.
 
We prove the connectedness of $\Stab \scD$ 
in \S \ref{section:connectedness_D}
and the injectivity of
$\rho$
in \S \ref{section:char_two}
and \S \ref{section:lifting}.
These results,
together with the above theorem of Bridgeland,
gives the following explicit description
of $\Stab \scD$:
%, at least under the assumption 
%that the base field is $\bC$:

\begin{theorem}
$\Stab \scD$ is the universal cover
of $\widehat{\frakh}^\reg / \widehat{W}$.
\end{theorem}

As for $\Stab \scC$,
a result of Thomas \cite{Thomas_SCBG}
shows that there is a distinguished connected component
of $\Stab \scC$
which is the universal cover of
$\frakh^\reg / W$.
We will prove the connectedness of $\Stab \scC$
in \S \ref{section:connectedness_C},
so that this connected component
is the whole of $\Stab \scC$.

%%%%%%%%%%%%%%%%%%%%%%%%%%%%%%%%%%%%%%%%%%%%%%%%%%%%%%%%%%%%%%%%%%%%%%%%
%%%%%%%%%%%%%%%%%%%%%%%%%%%%%%%%%%%%%%%%%%%%%%%%%%%%%%%%%%%%%%%%%%%%%%%%
\section{The McKay correspondence} \label{section:McKay}

We collect basic facts on the McKay correspondence
in this section.
We expect that the result in this section
is well-known to experts,
although we have been unable to locate
an appropriate reference.
Throughout this section,
$k$ will denote a field of any characteristic.
%unless otherwise specified.
We restrict our discussion
to the case of $A_n$-singularities
since it is the only case in need in this paper.
For a noetherian $k$-algebra $A$,
the abelian category of finitely generated
right $A$-modules
will be denoted by $\module A$.

\subsection{Path algebra and the endomorphism algebra of a reflexive 
module}
Let $\palg$
be the preprojective algebra
for the affine Dynkin quiver
of type $A_{n}^{(1)}$,
described explicitly
as follows:
As a $k$-vector space,
$\palg$ is generated by
the symbols
$(i_1|i_2|\dots|i_l)$
for $l\ge1$, $i_m \in 
\bZ/(n+1)\bZ$ and $i_{m+1}=i_m \pm 1$.
The multiplication is defined by
$$
(i_1|\dots|i_l)(j_1|\dots|j_m)=
\begin{cases}
(i_1|\dots|i_l|j_2|\dots|j_m) &\text{if } i_l=j_1, \\
0 &\text{otherwise,}
\end{cases}
$$
and the relations are generated by
$$
(i|i+1|i)=(i|i-1|i)
$$
for $i \in \bZ/(n+1)\bZ$.

Let $\scO=k[x, y, z]/(xy+z^{n+1})$ be the affine coordinate ring of the 
rational double point of type $A_{n}$.
For an integer $a=(n+1)q+r$ with $0\le r \le n$,
consider the fractional ideal $I_a =(y^{q+1}, y^qz^r)\scO$ of $\scO$.
$I_a$'s are reflexive $\scO$-modules such that
$I_a \cong I_b$ if and only if $a-b$ is divisible by $n+1$.
For $i \in \bZ/(n+1)\bZ$, we lift $i$ to $a \in \bZ$ with $0\le a \le n$ and 
put $E_i = I_a$.
For an integer $b=q(n+1)+r$ with $0\le r \le n$, we fix the isomorphism
$I_b \cong E_{(b\mod n+1)}=I_r$ given by the multiplication by $y^{-q}$.
Consider the reflexive $\scO$-module
$$
E= \bigoplus_{i \in \bZ/(n+1)\bZ} E_i.
$$
$E$ is an $\palg\otimes_k \scO$-module in the following way.
The idempotent $(i)$ acts as the projection of $E$ to $I_i$.
The path $(i|i+1)$ corresponds to the homomorphism $I_a \overset{\cdot 
z}{\longrightarrow} I_{a+1}$ given by the multiplication by $z$, where 
$a\in \bZ$ is a lift of $i$.
The path $(i|i-1)$ goes to the inclusion $I_a \hookrightarrow I_{a-1}$.
Then it is easy to see that we obtain an $\palg \otimes_k \scO$-action 
on $E$.
Thus we have a $k$-algebra homomorphism
$$
\eta: \palg \to \End_{\scO}(E).
$$
\begin{proposition}
$\eta$ is an isomorphism.
\end{proposition}
\begin{proof}
We first note that an element of $\palg$ is a $k$-linear combination of
the elements $P(i,l,m)$ defined as follows for $i \in \bZ/(n+1)\bZ$ and $l, 
m\in \bZ$ with $m \ge 0$:
\begin{equation*}
P(i,l,m)= \begin{cases}
(i|i+1|i)^m(i|i+1|\dots|i+l-1|i+l) & \text{if} \ l \ge 0, \\
(i|i+1|i)^m(i|i-1|\dots|i+l+1|i+l) & \text{otherwise}.
\end{cases}
\end{equation*}

To show that $\eta$ is an isomorphism, it is sufficient to show that
the restricted map
\begin{equation*}
\eta_{ij}: (i)\palg(j) \to \Hom_{\scO}(E_i, E_j)
\end{equation*}
is bijective for $i, j \in \bZ/(n+1)\bZ$.
Lift $i, j$ to $a, a+r \in \bZ$ with $0 \le r \le n$.
Then we have isomorphisms
$$
\Hom_{\scO}(E_i, E_j) \cong \Hom_{\scO}(I_a, I_{a+r}) \cong I_r=(y, 
z^r)\scO.
$$

If $l= q(n+1)+r$ with $q \ge 0$, then $P(i, l, m)$ is mapped to $x^qz^{m+r} 
\in I_r$,
and if $l = -q(n+1) +r<0$ with $q>0$, then $P(i, l, m)$ is mapped to 
$y^qz^{m} \in I_r$.
Moreover, the monomials $x^qz^{m+r}$ ($q \ge0, m \ge0)$ and $y^qz^{m}$ 
($q>0, m \ge0$)
form a $k$-linear basis of $I_r$.
Therefore, $\eta_{ij}$ must be bijective.
\end{proof}

\subsection{A full sheaf as a projective generator}
Let $f: X \to Y = \Spec \scO$ be the minimal resolution
and $\Per{-1}(X/Y)$ be the abelian category
of perverse sheaves
introduced by Bridgeland \cite{Bridgeland_FDC};
an object
$E \in \Per{-1}(X/Y)$
is a bounded complex of coherent sheaves on $X$
such that
its cohomology sheaves satisfy
$$ 
f_*(\scH^{-1}(E)) = 0,\quad R^1 f_*(\scH^{0}(E)) = 0,\quad
\scH^{i}(E) = 0 \mbox{ for } i \neq -1,0 ,
$$
and
$$
 \Hom_X(\scH^{0}(E), F) = 0
$$
for any coherent sheaf $F$ on $X$
satisfying $\bR f_* F = 0$.

For a reflexive $\scO$-module $F$,
put
$$
\wt{F}:=f^{*}(F)/\text{torsion},
$$
which is a locally free sheaf on $X$(see \cite{Artin-Verdier}).
A locally free sheaf of this form is called a {\it full sheaf}.
It is proved in \cite{Esnault} that a locally free sheaf $\scF$ on $X$
is a full sheaf if and only if the following two conditions are 
satisfied:
\begin{enumerate}
\item $\scF$ is generated by its global sections.
\item $R^1f_*(\scF^{\vee})=0$.
\end{enumerate}

The following result
is due to Van den Bergh:

\begin{proposition}
[{\cite[proposition 3.2.7, corollary 3.2.8]{Van_den_Bergh_TFNR}}]
\begin{enumerate}
\item
A full sheaf $\scM$
is a projective generator of $\Per{-1}(X/Y)$
if and only if
its first Chern class $c_1(\scM)$
is ample
and
$\scO_X$ is a direct summand
of $\scM^{\oplus a}$
for some positive integer $a$.
%\end{proposition}
%\begin{proposition}
%[{\cite[proposition 3.2.8]{Van_den_Bergh_TFNR}}]
\item
Assume that a full sheaf $\scM$
is a projective generator of $\Per{-1}(X/Y)$
and
put $A = \End_X(\scM)$.
Then the functor
$
 \RHom_X(\scM, \bullet)
$
gives an equivalence
between $D^b \coh X$
and $D^b \module A$,
whose inverse is given by
the functor
$
 \bullet \Lotimes_A \scM.
$
\end{enumerate}
\end{proposition}

The above proposition
yields the following:

\begin{theorem} \label{th:McKay}
The bounded derived category $D^b \coh X$
of coherent sheaves on the crepant resolution
$X$ of the $A_n$-singularity
is equivalent
to the bounded derived category
$D^b \module \palg$
of finitely generated
right $\palg$-modules.
\end{theorem}

\begin{proof}
Put $E=\oplus_{i \in \bZ/n\bZ} E_i$ be
as in the previous subsection.
Since
$$
 c_1(\wt{E_i}) \cdot C_j = \delta_{ij},
$$
the corresponding full sheaf $\wt{E}$
has an ample first Chern class
and is hence
a projective generator of $\Per{-1}(X/Y)$
satisfying $\End(\wt E) \cong \palg$.
\end{proof}

Let $\module_0 \palg$ be the abelian category
of finitely generated nilpotent right $\palg$-modules.
Under the above equivalence,
$\scD_k$ corresponds
to the bounded derived category $D^b \module_0 \palg$
of $\module_0 \palg$.

When $k$ contains a primitive $(n+1)$th root $\zeta$ of unity,
$Y$ is isomorphic to the quotient
$\bA^2 / G$ of the affine plane
by the natural action of the subgroup $G$
of $\SL_2(k)$ generated by the diagonal matrix
$\diag(\zeta, \zeta^n)$.
Since $\palg$ is isomorphic to the crossed product
$
 k[x, y] \rtimes k[G]
$
of the polynomial ring with the group ring,
the category of finitely generated nilpotent $\palg$-modules
is equivalent to the category $\coh_0^G \bA^2$
of $G$-equivariant coherent sheaves on $\bA^2$
supported at the origin:
$$
 \module_0 \palg \cong \coh_0^G \bA^2.
$$
Hence Theorem \ref{th:McKay} in this case
gives the equivalence
$$
 D^b \coh_Z X \cong D^b \coh_0^G \bA^2
$$
of triangulated categories,
first proved by Kapranov and Vasserot
\cite{Kapranov-Vasserot}
(see also Bridgeland, King, and Reid
\cite{Bridgeland-King-Reid}).

%%%%%%%%%%%%%%%%%%%%%%%%%%%%%%%%%%%%%%%%%%%%%%%%%%%%%%%%%%%%%%%%%%%%%%%%
%%%%%%%%%%%%%%%%%%%%%%%%%%%%%%%%%%%%%%%%%%%%%%%%%%%%%%%%%%%%%%%%%%%%%%%%
\section{Connectedness of $\Stab \scD$}
 \label{section:connectedness_D}

We prove the connectedness of $\Stab \scD$
in this section.
Our strategy is the following:
\begin{enumerate}
\item \label{item:standard_stability_I}
 Put
 \[
  U:=\bigl\{ \sigma\in \Stab \scD
  \bigm| \mbox{$\scO_x$ is $\sigma$-stable for all the closed points
   $x\in Z$
  }
  \bigr\}.
 \]
 %(A stability condition belonging to $U$ is sometimes called a
 %\emph{standard stability condition}.)
 Then $U$ is connected.
 The proof is parallel to
 the $K3$ case
 due to Bridgeland \cite[\S 11]{Bridgeland_SCK3}.
\item \label{item:standard_stability_II}
 The connected component of $\Stab \scD$
 containing $U$
 is preserved by the action of $\Br(\scD)$.
 This follows from \cite[theorem 1.1]{Bridgeland_SCKS}
 since $U$ belongs to the distinguished connected component
 studied by Bridgeland.
\item \label{item:perturbation}
 For any connected component $\Sigma$ of $\Stab \scD$,
 there is a stability condition $\sigma=(Z, \scP) \in \Sigma$ such that
 $Z(\scO_x) \in \bR_{< 0}$ for a closed point $x\in Z$
 and
 $Z(E) \not\in \bR$ for any spherical object $E$.
\item \label{item:technical_lemma}
 For a stability condition $\sigma$ as above,
 there is a $\sigma$-stable object
 $\omega \in \scP(1)$
 satisfying some technical conditions.
\item \label{item:finding_stable_point}
 The technical conditions together with
 a result in \cite{Ishii-Uehara_ADC}
 imply the existence
 of a closed point $x \in Z$,
 an autoequivalence $\Phi \in \Br(\scD)$,
 and an integer $d \in \bZ$
 such that
 $\omega = \Phi(\scO_x[d])$.
\item \label{item:induction}
 If $\scO_x$ is stable with respect to a stability 
 condition $\sigma$ for some point $x \in C_i$,
 then $\sigma$ induces stability conditions
 on $D^b \coh_{Z'} X$ and $D^b \coh_{Z''} X$,
 where $D^b \coh_{Z'} X$ and $D^b \coh_{Z''} X$
 are the full triangulated subcategories of $\scD$
 consisting of objects
 supported at $Z' = C_1 \cup \dots \cup C_{i-1}$
 and
 $Z'' = C_{i+1} \cup \dots \cup C_{n}$,
 respectively.
\item \label{item:connectedness}
 This proves the existence of
 an autoequivalence
 $\Phi \in \Br(\scD)$
 such that
 $\Phi \sigma$ belongs to $U$
 by induction on $n$.
\end{enumerate}
Note that
the above argument
except for (v)
works not only for $A_n$-singularities
but also for any rational double points
in dimension two.

%%%%%%%%%%%%%%%%%%%%%%%%%%%%%%%%%%%%%%%%%%%%%%%%%%%%%%%%%%%
\subsection{The connected component containing $U$}

Here we collect basic results
%of Bridgeland
on the connected component of $\Stab \scD$
containing a stability condition
where the structure sheaves of all the closed points
are stable.
The following lemma is taken from
\cite[lemma 10.1]{Bridgeland_SCK3}:

%%
%% lemma:K3

\begin{lemma}\label{lemma:K3}
For a stability condition $\sigma =(Z,\scP)$
and a closed point $x \in Z$,
%\red{with resudie field $k$}
assume that $\scO _x$ is $\sigma$-stable of phase $1$.
\begin{enumerate}
\item
If $E \in \scP((0,1])$, then $E$ is a sheaf in a neighborhood of $x$.
\item
If $E \in \scP(1)$ and $E$ is stable,
then $E=\scO_x$ or $x \not\in\Supp E$.
\end{enumerate}
\end{lemma}

The following two lemmas are also essentially
due to Bridgeland,
whose proofs we include for completeness.
The proof of Lemma \ref{lemma:stability_i}
is similar to,
but easier than,
the corresponding statement for K3 surfaces
in \cite[\S 11]{Bridgeland_SCK3}.
As for Lemma \ref{lemma:stability_ii},
we use the result
in \cite{Bridgeland_SCKS}
rather than imitate
the proof in \cite[\S 12]{Bridgeland_SCK3}.

\begin{lemma} \label{lemma:stability_i}
 The subset $U \subset \Stab \scD$ defined by
 \[
  U:=\bigl\{ \sigma\in \Stab \scD
  \bigm| \mbox{$\scO_x$ is $\sigma$-stable
   for all the closed points $x\in Z$
  }
  \bigr\}
 \]
 is connected.
\end{lemma}

\begin{proof}
Let us fix a point $y \in Z$.
Consider a stability condition $(Z, \phi) \in U$ that satisfies $\scO_y \in \scP(1)$ and $Z(\scO_y)=-1$.
Then, for any point $x \in Z$, $Z(\scO_x)=Z(\scO_y)=-1$ implies $\phi(\scO_x)\in 2\bZ+1$.
Suppose $y$ is contained in an irreducible component $C_i$ of $Z$.
Then the support of at least one Harder-Narasimhan factor of $\scO_{C_i}$ contains $C_i$.
Hence there is a stable object $E \in \scP((0,1])$ whose support contains $C_i$.
$E$ is a sheaf in a neighborhood of $y$ by Lemma \ref{lemma:K3} and therefore
for a point $x$ in a neighborhood of $y$ in $C_i$, $\phi(\scO_x)$ must be $1$. 
Since $Z$ is connected, $\phi(\scO_x)=1$ holds for every point $x \in Z$
and $(Z, \phi)$ belongs to
the following subset $V$ of $U$:
$$
 V := \left\{ \left. (Z, \scP) \in U \, \right| \,
    \scO_x \in \scP(1) \ \text{and} \ 
     Z(\scO_x) = - 1
    \ \text{for any point} \ x \in Z
      \right\}. 
$$
Therefore
we can take any stability conditions in $U$
into $V$
by the action of the connected group
$\GLt(2, \bR)$.
Hence
it suffices to show that
$V$ is connected.
Take a stability condition $\sigma=(Z,\scP)\in V$.
Lemma \ref{lemma:K3} implies that the heart $\scP((0,1])$
coincides with the abelian category $\coh_Z X$
and every one-dimensional sheaf in $\coh_Z X$
has the phase $\phi$ with $0<\phi<1$.
It follows that
$Z$ belongs to
the subset $\scL$ of
$\Hom(K(\scD),\bC)$ 
defined by
\[
\scL:=
\bigl\{ Z_{(\beta, \omega)} \in \Hom(K(\scD),\bC)
\bigm| \beta\in N^1(X/Y), \omega\in\scA(X/Y) \bigr\}.
\]
Here
$\scA(X/Y)$ is the ample cone in $N^1(X/Y)$ and
$$
 Z_{(\beta, \omega)}(E)
  = - ch_2(E) + (\beta + \sqrt{-1} \omega) \cdot c_1(E)
$$
for $E \in K(\scD)$.
Now consider the restriction $\scZ |_V:V\to \scL$
of the map $\scZ$ in Theorem \ref{theorem:annals},
which is clearly injective.
We will show in the proof of 
Lemma \ref{lemma:perturb} that 
$V(\Sigma)=\Hom(K(\scD), \bC)$  
for any connected component $\Sigma$ of $\Stab \scD$.
Therefore $\scZ |_V$ is also a local homeomorphism.
To show that it is surjective,
take any $Z_{(\beta, \omega)}$ in $\scL$.
Since $\omega$ is ample,
$Z_{(\beta, \omega)}$ is a stability
function on $\coh_Z X$.
To show that it has the Harder-Narasimhan property,
we check the conditions (a) and (b)
in \cite[proposition 2.4]{Bridgeland_SCKS}.
Condition (a) follows from the fact that
for any infinite sequence
$$
 \dots \subset E_{j+1} \subset E_j \subset \dots
  \subset E_2 \subset E_1
$$
of proper subobjects,
there is a natural number $N$ such that
for any $l > m > N$,
the difference $[E_m] - [E_l]$ is
a positive multiple of $\delta$
defined in page \pageref{def:delta}.
Since $\coh_Z X$ is noetherian,
there is no infinite sequence of proper quotients
and the condition (b) is automatic.
Local-finiteness also follows from a similar reasoning.
Hence $\scZ|_V$ gives a homeomorphism
from $V$ to $\scL$,
which is connected.
\end{proof}

\begin{lemma} \label{lemma:stability_ii}
The connected component of $\Stab \scD$
containing $U$
is preserved
by the action of $\Br(\scD)$.
\end{lemma}

\begin{proof}
Take $\omega \in \scA(X/Y)$ satisfying
$\omega \cdot C_j =1$ for all $j$.
Then there is a unique stability condition
$\sigma = (Z, \scP) \in V$
such that
$Z= \exp(\sqrt{-1} \omega)$.
This stability function satisfies
$Z(\scO_x)=-1$,
$Z(\scO_{C_j}(-1))=\sqrt{-1}$ for all $j$
and $Z(\omega_Z)=-1+n \sqrt{-1}$.
Moreover,
we can easily check that
the sheaves $\scO_{C_j}(-1)$ for $j=1, \dots, n$
and $\omega_Z$
are $\sigma$-stable.
Thus, if we take $\alpha \in (0, 1/2)$ with
$\tan (\pi \alpha) > n$,
then the abelian category
$\scP((\alpha, \alpha+1])$ contains
$\omega_Z[1], \scO_{C_1}(-1), \dots, \scO_{C_n}(-1)$,
which corresponds to
the simple $\palg$-modules
by the McKay correspondence.
Since both are the hearts of bounded t-structures,
they must coincide:
$$
 \scP((\alpha, \alpha+1])=\module \palg.
$$
This shows that
$\sigma$ lies in the distinguished connected component
appearing in the work of Bridgeland \cite{Bridgeland_SCKS}.
Hence it is preserved by the action of $\Br(\scD)$
\cite[theorem 1.4]{Bridgeland_SCKS}.
\end{proof}

%%%%%%%%%%%%%%%%%%%%%%%%%%%%%%%%%%%%%%%%%%%%%%%%%%%%%%%%%%%
\subsection{Perturbation lemma}

Here we prove the following lemma.
Recall that $\delta \in K(\scD)$ is
the class of the structure sheaf
of a closed point. %with residue field $k$.
%in page \pageref{def:delta}.

\begin{lemma}\label{lemma:perturb}
For any connected component $\Sigma$ of $\Stab \scD$,
there is a stability condition $\sigma=(Z, \scP) \in \Sigma$ such that
$Z(\delta) \in \bR_{< 0}$
and
$Z(E) \not\in \bR$ for any spherical objects $E$.
\end{lemma}

\begin{proof}
Recall that
in the proof of Theorem \ref{theorem:annals},
Bridgeland describes $V(\Sigma)$
explicitly as
$$
 V(\Sigma)=\{U \in \Hom(K(\scD), \bC)\,|\, \|U\|_{\sigma} < \infty
%   \ \text{for any} \ \sigma \in \Sigma
           \},
$$
where
$$
\|U\|_{\sigma}= \sup \left\{\left. \frac{|U(E)|}{|Z(E)|} \right|
\text{$E$ is stable in $\sigma$} \right\}
$$
for a stability condition
$\sigma=(Z, \scP) \in \Sigma$.
He also shows that
$V(\Sigma)$ does not depend
on the choice of $\sigma$.

%%
%% {lemma:vsigma}

Now we show that
for any connected component
$\Sigma$ of $\Stab \scD$,
$V(\Sigma)$ coincides with $\Hom(K(\scD), \bC)$.
Let $\sigma=(Z, \scP)$ be a stability condition in $\Sigma$.
We want to show that $\|U\|_{\sigma} <\infty$
for every $U \in \Hom(K(\scD), \bC)$.

We first consider the case where $Z(\delta)=0$.
%(Actually, it will turn out that this is not the case,
%but we don't prove it here.)
In this case, every stable object must be spherical
and its class in the finite root lattice $K(\scD)/\bZ \delta$
is a root.
Assume that $E_1$ and $E_2$ are stable objects
from the heart $\scP((0,1])$ such that
their classes in $K(\scD)/\bZ \delta$ coincide.
Then we have $\chi(E_1, E_2)=2$ and $Z(E_1)=Z(E_2)$.
It follows that $E_1$ and $E_2$ have the same phase and
there is a non-zero morphism between them.
If this happens for stable objects $E_1$ and $E_2$,
then we must have $E_1 \cong E_2$.
Thus for every root in $K(\scD)/\bZ \delta$,
there is at most one corresponding stable object in $\scP((0,1])$.
This shows that there are only finitely many isomorphism classes
of stable objects in $\scP((0,1])$.
Therefore, we have $\|U\|_{\sigma} <\infty$
for every $U$.

Consider the remaining case: $Z(\delta) \ne 0$.
For a stable object, its class in the affine root lattice $K(\scD)$ is a root.
Hence there are finitely many stable objects $E_1, \dots, E_m$ such that
for every stable $E$, there exists $E_j$ with
$[E]-[E_j] \in \bZ \delta \subseteq K(\scD)$.
Since
\[\lim_{k \to \pm \infty} \frac{U([E_j] + k \delta)}{Z([E_j] + k \delta)}
=\frac{U(\delta)}{Z(\delta)}\]
converges for any $U \in \Hom(K(\scD), \bC)$,
we obtain $\|U\|_{\sigma} <\infty$.

Therefore,
we can find a stability condition
$\sigma' = (Z', \scP') \in \Sigma$
such that
$Z'(\delta) \neq 0$.
Then we can bring $\sigma'$
to
$\sigma'' = (Z'', \scP'') \in \Sigma$
such that
$Z''(\delta) \in \bR_{<0}$
by the action
of $\GLt(2, \bR)$.
Since
the set of classes of spherical objects
in $K(\scD) / \bZ \delta$ is finite,
we can perturb $\sigma''$ further
to
$\sigma = (Z, \scP)$
satisfying
$Z(E) \not\in \bR$
for any spherical objects $E$,
while keeping the condition
$Z(\delta) \in \bR_{<0}$ intact.
\end{proof}

If $\sigma$ satisfies the conditions
in Lemma \ref{lemma:perturb},
then the imaginary part
$\Im Z$
defines a Weyl chamber
of the finite root lattice  
$K(\scD)/\bZ \delta$;
the
simple roots are
linearly independent elements
$\alpha_1, \dots, \alpha_n \in  
K(\scD)/\bZ \delta$
with $\Im Z(\alpha_i) >0$ such that
any root $\alpha \in K(\scD)/\bZ \delta$
with $\Im Z(\alpha) >0$ is a linear combination of $\alpha_1, \dots,  
\alpha_n$
whose coefficients are non-negative.

%%%%%%%%%%%%%%%%%%%%%%%%%%%%%%%%%%%%%%%%%%%%%%%%%%%%%%%%%%%
\subsection{Main technical lemma}

Here we prove Lemma \ref{lemma:rose},
which lies at the heart of our proof
of the connectedness of $\Stab \scD$.

%%
%% lemma:rose

\begin{lemma} \label{lemma:rose}
Assume that a stability condition $\sigma=(Z, \scP)$ satisfies
the conditions
in Lemma \ref{lemma:perturb}.
Then there are $\sigma$-stable spherical objects $E, E'\in \scP((0,1))$ with  
$\Hom^1_{\scD}(E,E')=0$
and a $\sigma$-stable object $\omega \in \scP(1)$
which fit into the following
short exact sequence in $\scP((0,1])$;
\begin{equation}\label{eq:omega}
0\to E'\to E\to \omega \to 0.
\end{equation}
\end{lemma}

\begin{proof}
The proof goes in five steps:

%%
%% {lemma:toroidal}

\begin{step} \label{step:stable_toroidal}
There is a
$\sigma$-stable object $\omega \in \scP(1)$
such that $[\omega] \in \bZ \delta$.
\end{step}

We will show that
there is a non-zero object in $\scP(1)$,
and the desired stable object in $\scP(1)$
will be obtained
as its Jordan-H\"{o}lder component.

First note that
there are objects $E_1, E_2$
of the heart
$\scA := \scP((0,1])$
%of the bounded $t$-structure
%corresponding to $\sigma$
%by \cite[proposition 5.3]{Bridgeland_SCTC}
such that
$[E_1] - [E_2] = \delta$
since $K(\scA) \cong K(\scD)$.
%If $E_1 = 0$,
%then $0 \ne E_2 \in \bZ \delta$
%and hence $E_2 \in \scP(1)$.
%Assume that $E_1 \ne 0$.
If $E_1$ is in $\scP(1)$,
then we are done.
If $E_1 \not \in \scP(1)$,
then $[E_1] \not \in \bZ \delta$
and hence we have
$$
 \chi(E_1, E_2)
  = \chi(E_1, E_1 - \delta)
  = \chi(E_1, E_1)
 \ge 2.
$$
It follows that at least one of
$\Hom_{\scD}(E_1, E_2)$
and
$\Hom_{\scD}(E_2, E_1)$
is non-zero.
First assume that
there is a non-zero morphism
$
 f : E_1 \to E_2
$
and
let $F$ be the image of $f$
in the abelian category $\scP((0,1])$.
If $F \in \scP(1)$, then we are done.
If not, put
$
 E_1' := \Ker f$
and
$
 E_2':=\Coker f
$
in $\scP((0,1])$.
Then
$
 [E_1'] - [E_2'] = [E_1] - [E_2] \in \bZ \delta
$
and we can repeat the argument.
Since $[E_1] - [E_1'] = [F]$
and the images of $[E_1], [E_1']$, and $[F]$
in the finite root lattice
$K(\scD)/\bZ \delta$
are sums of positive roots,
this process terminates
after finitely many steps.
The case 
$\Hom_{\scD}(E_2, E_1)\ne 0$
can be treated similarly.
Hence there is a stable object in $\scP(1)$.
Since we assume $\sigma$ satisfies the conditions in Lemma \ref{lemma:perturb},
the class of an object in $\scP(1)$ lies in $\bZ \delta$ and Step 1 is proved.

We impose an additional condition on $\omega$,
which will be crucial in Step 5.
Put
$$
 a = \mathrm{min} \{ l \in \bZ \mid 
 \text{ there is a non-zero object }\omega \in \scP(1)
      \text{ such that } [\omega] = l \delta \}
$$
and
let $\omega \in \scP(1)$ be an object
such that $[\omega] = a \delta$.
The miminality of $a$ implies that $\omega$ is a simple object of $\scP(1)$,
which means $\omega$ is stable.

\begin{step} \label{step:orthogonal_decomposition}
There is an object $H \in \scP((0,1))$
such that
$\Hom_{\scD}(H, \omega) \ne 0$.
\end{step}

Assume that
$
 \Hom_{\scD}(H, \omega)=0
$
for any $H \in \scP((0,1))$.
Since $\omega \in \scP(1)$
and $H \in \scP((0,1))$,
we also have
$$
 \Hom_{\scD}^{\leq 0}(\omega, H)
  = \Hom_{\scD}^{\leq -1}(H, \omega)
  = 0.
$$
Moreover, it follows from $[\omega] \in \bZ \delta$ that
$\chi(\omega, H)=0$.
Then the Serre duality shows that
$$
 \Hom^*_{\scD}(H, \omega) = \Hom^*_{\scD}(\omega, H) = 0.
$$
For a stable object
$\omega' \in \scP(1)$
which is not isomorphic to $\omega$,
we have
$$
 \Hom_{\scD}^{\leq 0}(\omega', \omega)
 = \Hom_{\scD}^{\leq 0}(\omega,  \omega')
 = 0.
$$
This, together with
$
 \chi(\omega, \omega')=0
$
and the Serre duality,
again implies
$$
 \Hom^*_{\scD}(\omega', \omega)
  = \Hom^*_{\scD}(\omega, \omega')
  = 0.
$$
Hence $\scD$ admits an orthogonal decomposition
$$
 \scD = \langle \Omega, \Omega^{\perp} \rangle
$$
into $\Omega = \langle \omega \rangle$
and $\Omega^\perp = \! \ ^{\perp} \Omega$.
This is impossible
since $Z$ is connected
(cf. \cite[example 3.2]{Bridgeland_FMTES}).

\begin{step} \label{step:any_quotient}
There is an object $E \in \scP((0,1))$ such that
$\Hom_{\scD}(E, \omega) \ne 0$
and
for any non-zero subobject $F \subset E$
in the abelian category $\scP((0,1])$,
the quotient object $E/F$ lies in $\scP(1)$.
\end{step}

If $H$ in Step \ref{step:orthogonal_decomposition}
satisfies the second condition
in Step \ref{step:any_quotient},
put $E = H$.
If not, then there is a non-zero subobject $F \subset H$
such that $H/F$ does not lie in $\scP(1)$.
If $H/F$ is in $\scP((0,1))$, then put $H_1=F$
and $H_2=H/F$.
If $H/F$ is not in $\scP((0,1))$,
then let $H_1 \subset H$ be the pull-back (to $H$)
of the first Harder-Narasimhan factor of $H/F$
and put $H_2=H/H_1$.
Both $H_1$ and $H_2$ are in $\scP((0,1))$ and
either $\Hom_{\scD}(H_1, \omega)\ne 0$
or
$\Hom_{\scD}(H_2, \omega) \ne 0$
holds.
If $\Hom_{\scD}(H_i, \omega) \ne 0$, we can replace $H$ by $H_i$.
We have $[H]=[H_1]+[H_2]$ in $K(\scD)$ and the images of these objects
in $K(\scD)/\bZ \delta$ are non-zero sums of positive roots.
Thus we cannot repeat the process infinitely many times and
Step \ref{step:any_quotient} is proved.

\begin{step}
Let $f: E \to \omega$ be a non-zero morphism
and put $E':= \Ker f$ in $\scP((0,1])$.
Then $E$ and $E'$ are stable and spherical.
\end{step}

Since $E$ and $E'$ are not in $\scP(1)$,
Lemma \ref{lemma:stable_real_root_spherical} implies that
it suffices to show
that they are stable.
Let $0 \ne F \subsetneq E$ be a proper subobject.
Then by Step \ref{step:any_quotient}, $E/F$ is in $\scP(1)$.
From the equality
$$
 Z(F)=Z(E)-Z(E/F)
$$
and the fact that
$$
 \Im Z(F) > 0, \ 
 \Im Z(E) > 0, \ 
 \text{and} \ 
 Z(E/F) \in \bR_{<0},
$$
we obtain $\phi(F) < \phi(E)$.
The proof of the stability of $E'$
is the same.

\begin{step}
$\Hom^1_{\scD}(E, E')=0$.
\end{step}

Assume $\Hom^1_{\scD}(E, E') \ne 0$.
Then there is a non-trivial extension
\begin{equation}\label{equation:stableextension}
0 \to E' \to L \to E \to 0
\end{equation}
in $\scP((0,1])$.
We will show that
$L$ is stable,
which contradicts
$\chi(L, L)=8$
and Lemma \ref{lemma:root}.
Let $M \subsetneq L$
be a non-zero proper subobject
and
put $F':=E' \cap M \subset E'$
and $F:=M/F' \subset E$;
$$
\begin{array}{ccccccccc}
 0 & \to & F' & \to & M & \to & F & \to & 0 \\
 & & \cap
 & & \cap
 & & \cap
 & & \\
 0 & \to & E' & \to & L & \to & E & \to & 0. \\
\end{array}
$$
Then there are three cases to be considered:
\begin{itemize}
 \item {\it The case $F \ne 0$ and $F' \ne 0$}:

By Step \ref{step:any_quotient},
we have $E'/F', E/F \in \scP(1)$.
From the equality
$$
 Z(M) = Z(L) - (Z(E'/F')+Z(E/F))
$$
and the fact that
$Z(E'/F'), Z(E/F) \in \bR_{<0}$,
we obtain $\phi(M) < \phi(L)$.

\item {\it The case $F=0$ and $F' \ne 0$}:

Since $E'$ is stable and $Z(L)=2Z(E')+Z(\omega)$,
we have $\phi(M)=\phi(F') \le \phi(E') < \phi(L)$.

\item {\it The case $F \ne 0$ and $F' = 0$}:

If $F=E$, then \eqref{equation:stableextension}
splits and contradicts our assumption.
Hence we have $F \ne E$
and
$[E / F]=b \delta$
with $b > 0$
by Step \ref{step:any_quotient}.
Then our choice of $\omega$ implies $a \le b$.
From the equality
$$
 Z(M)
  = Z(F)
  = Z(E) - b Z(\delta)
  = Z(E') + (a - b) Z(\delta),
$$
we conclude $\phi(M) \le \phi(E') < \phi(L)$.

\end{itemize}

\end{proof}

%%%%%%%%%%%%%%%%%%%%%%%%%%%%%%%%%%%%%%%%%%%%%%%%%%%%%%%%%%%
\subsection{Stabilizing a skyscraper sheaf}

The proof of \cite[proposition 1.7]{Ishii-Uehara_ADC}
actually shows the following  
stronger result. 
%Note that the proof there works over any field $k$. 

%%
%% {lemma:IU05}

\begin{lemma} \label{lemma:IU05}
Let $E$ and $E'$ be
spherical objects in $\scD$
satisfying
$$
 \Hom^{i}_{\scD}(E',E)
  \cong \begin{cases}
          \bC^{\oplus 2} & \text{if } i = 0, \\
            0 & \text{otherwise}.
        \end{cases}
$$
Then there is an autoequivalence
$\Phi \in \Br(\scD) \times \bZ$
such that
$\Phi(E')\cong \scO_{C_l}(a-1)$
and
$\Phi(E)\cong \scO_{C_l}(a)$
for some $l \in \{ 1, \dots, n \}$
and $a \in \bZ$.
\end{lemma}

Let $E, E', \omega$ be the objects obtained in Lemma \ref{lemma:rose}.
Then we obtain an autoequivalence $\Phi \in \Br(\scD) \times \bZ$
by applying Lemma \ref{lemma:IU05} to $E$ and $E'$.
Since $\omega$ fits into the short exact sequence \eqref{eq:omega},
we have
$$
\Phi(\omega) \cong \{\scO_{C_l}(a-1) \to \scO_{C_l}(a)\}
\cong \scO_y
$$
for some point $y \in C_l$.
Replacing $\omega$ by its shift if necessary, we obtain the following:

%%
%% proposition:key

\begin{proposition}\label{proposition:key}
Let $\sigma = (Z, \scP)$ be
a stability condition
satisfying the conditions
in Lemma \ref{lemma:perturb}.
Then there is a $\sigma$-stable object
$\omega \in \scD$
such that
$\Phi(\omega) \cong \scO_y$
for a closed point $y \in Z$
%with residue field $k$
and an autoequivalence
$\Phi \in \Br(\scD)$.
\end{proposition}

%%%%%%%%%%%%%%%%%%%%%%%%%%%%%%%%%%%%%%%%%%%%%%%%%%%%%%%%%%%
\subsection{Induction lemma}

Here
we prove Corollary \ref{corollary:induction},
which states that
a stability condition $\sigma$
which stabilizes the skyscraper sheaf $\scO_x$
for some closed point $x \in C_i$
induces a stability condition
on the subcategory of $\scD$
consisting of objects
supported at $Z' = \bigcup_{i \ne j} C_j$.
The following lemma is a consequence of
Lemma \ref{lemma:K3}:

\begin{lemma}\label{lemma:march}
Assume that $\scO_x$ is $\sigma$-stable for some 
stability condition $\sigma$. For an object $E \in \scD$,
let
\begin{equation*}
\psmatrix[colsep=.5, rowsep=.7]
  0 = \Rnode{a}{E_0} & & E_1 & & E_2 & \cdots & E_{m-1} & &  
\Rnode{n}{E_m}=E \\
        & A_1 & & A_2 &   &        &       & A_m & \\
\endpsmatrix
\psset{nodesep=3pt,arrows=->}
\ncline{a}{1,3}
\ncline{1,3}{1,5}
\ncline{1,6}{1,7}
\ncline{1,7}{n}
\ncline{1,3}{2,2}
\ncline{1,5}{2,4}
\ncline{n}{2,8}
\psset{arrows=-}
\ncline{1,5}{1,6}
\psset{arrows=->,linestyle=dashed}
\ncline{2,2}{a}
\ncline{2,4}{1,3}
\ncline{2,8}{1,7}
\end{equation*}
be a collection of triangles
such that
for some real numbers
$$
 \phi_1 \geq \phi_2 \geq \dots \geq \phi_m ,
$$
$A_j$ is in $\scP(\phi_j)$ and stable
for $j = 1, \dots, m$,
obtained by combining the Harder-Narasimhan filtration
with a Jordan-H\"{o}lder filtration.
If $\Supp E $ does not contain $x$, then
$\bigcup_{j=1}^m \Supp A_j$ does not contain $x$ either.
\end{lemma}

\begin{proof}
We may assume that the phase of $\scO_x$ is $1$
and $E$ is not  $\sigma$-stable.
Suppose for contradiction that $x \in \Supp A_j$ for some $j$.
%Take the maximal $j$ in the set of the indices
%of such $A_j$'s, and replace $E_j$ with $E$.
By replacing $E$ with $E_j$ for a suitable $j$,
we may assume that $x \in \Supp A_m$.
Since the morphism $E \to A_m$ is not zero,
we have $\Supp A_m \supsetneq \{ x \}$.
Then Lemma \ref{lemma:K3}(ii) ensures that
$\phi(A_m)$ is not an integer and hence
$d < \phi(A_m) < d+1 $ for some $d \in \bZ$.
Since $A_m[-d]$ is a sheaf near $x$ by Lemma \ref{lemma:K3}(i)
and $\Supp A_m \ni x$ by our assumption,
we obtain $\Hom_{\scD}(A_m, \scO_x[d]) \ne 0$,
which implies
$\Hom_{\scD}(E_{m-1}, \scO_x[d-1]) \ne 0$.
This shows that
$\Hom_{\scD}(A_l, \scO_x[d-1]) \ne 0$
for some $l$ with $l < m$.
Hence we have
$\phi(A_l) \le \phi(\scO_x[d-1])=d$,
which contradicts the fact that
$\phi(A_l) \ge \phi(A_m) > d$.
\end{proof}

The above lemma yields the following
two corollaries:

%%
%%  {corollary:induction0}

\begin{corollary}\label{corollary:induction0}
%\begin{enumerate}
%\item
Suppose that $\scO_x$ is $\sigma$-stable for some $x \in C_i$.
If a point $y \in C_i$
is not contained in $C_j$ for any $j \ne i$,
then $\scO_y$ is also $\sigma$-stable.
%\item
%Assume that $\scO_x$ and $\scO_y$ are $\sigma$-stable for some
%  $x \in C_i$, $y \in C_{j}$ with $i \ne j$.
%Then $\scO_z$ is $\sigma$-stable for $z\in C_i\cap C_{j}$.
%\end{enumerate}
\end{corollary}

\begin{proof}
We may assume that $x\ne y$.
Consider the above filtration for $E=\scO_y$.
The support of each stable factor $A_j$ is connected by its stability
and does not contain $x$ by the above lemma.
It follows that $\Supp A_j$ is just the single point $y$
or it does not contain $y$.
On the other hand, since $\scO_y$ is indecomposable,
the union $\bigcup_j \Supp A_j$ is connected.
Therefore, we have $\bigcup_j \Supp A_j=\{y\}$.
This shows that $A_j$ must be of the form $\scO_y[d]$ for some $d\in\bZ$
since we have $\Hom_{\scD}(A_j, A_j) \cong \bC$
and
$\Hom^{\le -1}_{\scD}(A_j, A_j) = 0$
by the stability of $A_j$.
Hence $\scO_y\cong A_j[-d]$ is stable.
%The proof of (ii) is similar.
\end{proof}

%%
%% corollary:induction

\begin{corollary} \label{corollary:induction}
Assume that $\scO_x$ is $\sigma$-stable for some $x \in C_i$,
and put $Z' = \bigcup_{j \ne i}C_j$.
Then $\sigma$-semistable factors of an object
supported at $Z'$
are again supported at $Z'$.
\end{corollary}

\begin{proof}
Corollary \ref{corollary:induction0} implies that
for any $y \in Z \setminus Z'$, $\scO_y$ is $\sigma$-stable.
Then the statement follows from Lemma \ref{lemma:march}.
\end{proof}

%%%%%%%%%%%%%%%%%%%%%%%%%%%%%%%%%%%%%%%%%%%%%%%%%%%%%%%%%%%
\subsection{The proof of the connectedness of $\Stab \scD$}
\label{subsection:proof_of_the_connectedness}

Here we finish the proof of the connectedness of $\Stab \scD$.

\begin{lemma}\label{lemma:U}
For any connected component $\Sigma$ of $\Stab \scD$,
there is a stability condition $\sigma \in \Sigma$
and an autoequivalence $\Phi \in \Br(\scD)$
such that
$\Phi \sigma$ belongs to $U$.
\end{lemma}

\begin{proof}
We use the induction on $n$.
The perturbation lemma guarantees
the existence of a stability condition $\sigma=(Z, \scP) \in \Sigma$  
satisfying
$Z(\delta) \in \bR_{< 0}$ and $Z(E) \not\in \bR$
for any spherical object $E$,
and
we wish to find
an autoequivalence $\Phi\in \Br(\scD)$
such that
$\Phi\sigma$ belongs to $U$.

Proposition \ref{proposition:key} says that
the structure sheaf $\scO_y$ of some point $y\in Z$
is stable with respect to $\Phi\sigma$ for some $\Phi\in  
\Br(\scD)$.
In the case $n=1$, at this stage,
Corollary \ref{corollary:induction0}
ensures that $\Phi\sigma$ belongs to $U$.
In the case $n>1$, assume $y \in C_l$ and
put $Z_1=C_1\cup\cdots\cup C_{l-1}$ and  
$Z_2=C_{l+1}\cup\cdots\cup C_n$.
Moreover we define
$$\Br(\scD)_1:=\Span{ T_{\scO_{C_i}(a)} \bigm| 1\le i\le l-1, a\in\bZ  
}$$ and
$$\Br(\scD)_2:=\Span{ T_{\scO_{C_i}(a)} \bigm| l+1\le i\le n, a\in\bZ  
}.$$
Corollary \ref{corollary:induction} implies that $\Phi\sigma$ naturally
induces stability conditions on $D^b \coh_{Z_1} X$ and
on $D^b \coh_{Z_2} X$,
respectively.
Then by the induction hypothesis,
we can find autoequivalences $\Phi_1\in \Br(\scD)_1$ and
$\Phi_2\in \Br(\scD)_2$ such that
$\Phi_1\Phi_2 \Phi \sigma \in U$. Note that
we can regard $\Phi_k$ as an element of $\Auteq \scD$
via the natural inclusion
$\Br(\scD)_j \subset \Br(\scD)$
for $j = 1, 2$.
\end{proof}

Lemmas \ref{lemma:stability_i}, \ref{lemma:stability_ii},
and \ref{lemma:U} imply the connectedness of $\Stab \scD$.

%%%%%%%%%%%%%%%%%%%%%%%%%%%%%%%%%%%%%%%%%%%%%%%%%%%%%%%%%%%%%%%%%%%%%%%%
%%%%%%%%%%%%%%%%%%%%%%%%%%%%%%%%%%%%%%%%%%%%%%%%%%%%%%%%%%%%%%%%%%%%%%%%
\section{Faithfulness in characteristic two}
 \label{section:char_two}

%
%-------- for the mirror part -------------------
%

\newcommand{\DG}{\ensuremath{DG}}
\newcommand{\PreTr}{\ensuremath{\operatorname{Pre-Tr}}}
\newcommand{\Lag}{\mathcal{L}ag}
\newcommand{\Lagt}{\widetilde{\Lag}}
\newcommand{\CM}{\mathcal{CM}}
\newcommand{\Res}{\operatorname{Res}}
\newcommand{\Igr}{\mathop{I^{\mathrm{gr}}}\nolimits}
\newcommand{\ctilde}{\widetilde{c}}
\newcommand{\Ltilde}{\widetilde{L}}
\newcommand{\stilde}{\widetilde{s}}
\newcommand{\scFtilde}{\Fuk W^{\sim}}
\newcommand{\alphatilde}{\widetilde{\alpha}}

%
%----------- text starts ---------------------------------------
%

In this section,
we prove the faithfulness
of the action of 
the affine braid group $\ABG$
on $\scD_k$
in characteristic two.
The assumption on the characteristic is needed
in Theorem \ref{theorem:Khovanov-Seidel}
due to Khovanov and Seidel
\cite{Khovanov-Seidel}.
Throughout this section,
$k$ will be a field
of characteristic two.
We adopt the following approach,
which generalizes
to the affine case
the method of
Khovanov, Seidel, and Thomas
\cite{Khovanov-Seidel, Seidel-Thomas}
who proved the faithfulness
of the braid group actions:
\begin{enumerate}
\item
 Consider the affine manifold
 $$
  W = \{ (x, y, z) \in \bC^2 \times \bCx | x y z = z^{n+1} - 1 \}
 $$
 with a natural exact symplectic structure.
 It is a conic fibration
 over the $z$-plane
 whose discriminant set is the set of $(n+1)$th roots
 of unity.
 Any curve on the $z$-plane
 without self-intersections
 starting and ending
 at the discriminant set
 gives an exact Lagrangian two-sphere
 of $W$.
\item
 There are $n+1$ Lagrangian two-spheres
 $\{ L_i \}_{i=0}^n$
 of $W$
 constructed from straight line segments
 $\{ c_i \}_{i=0}^n$
 on the $z$-plane
 connecting the neighboring discriminant points.
 Let $\Fuk W$ be the Fukaya category
 whose set of objects is $\{ L_i \}_{i=0}^n$
 and whose spaces of morphisms are
 Lagrangian intersection Floer complexes.
 Although Fukaya categories are not honest categories
 but only $A_\infty$-categories in general,
 the above $\Fuk W$ turns out
 to be a differential graded category
 with the trivial differential,
 since there are no non-constant
 pseudoholomorphic maps
 from a disk
 to $W$
 with Lagrangian boundary conditions
 given by $\{ L_i \}_{i=0}^n$.
 \item Let $D^b \Fuk W$ be the derived category of $\Fuk W$
 defined using twisted complexes.
 It is equivalent to the derived category
 $D^b \module_0 \palg$
 of finite-dimensional nilpotent representations of
 the path algebra $\palg$ with relations
 appearing in \S \ref{section:McKay}.
 This is the {\em homological mirror symmetry}
 of Kontsevich
 \cite{Kontsevich_HAMS}
 for $A_n$-singularities.
\item
 There are two ways to define
 an action of the affine braid group
 on $D^b \Fuk W$;
 one is algebraic and given
 by the twist functor,
% appearing in the work of Seidel and Thomas \cite{Seidel-Thomas},
 and the other is geometric
 and given by the symplectic Dehn twist.
% which originates
% from the note of Arnold \cite{Arnold_SMMF}.
 They coincide by Seidel
 \cite[proposition 9.1]{Seidel_K3},
 which allows us to use symplectic geometry of $W$
 to prove the faithfulness
 of the affine braid group action.

\item \label{it:commutativity_of_ABG_actions}
 The affine braid group also acts
 on curves on the $z$-plane.
 The actions of the affine braid group
 on curves on the $z$-plane
 and Lagrangian submanifolds of $W$
 commute with
 the construction of Lagrangian submanifolds
 from curves.

\item \label{it:floer_gin}
 By Khovanov and Seidel
 \cite[theorem 1.3]{Khovanov-Seidel},
 the dimension of the Floer cohomology group
 between two Lagrangian two-spheres
 coming from two curves on the $z$-plane
 is given by twice their geometric intersection number.
 This step requires us to work in characteristic two.
\item \label{it:gin_faithful}
 Suppose given an element $b$ of the affine braid group
 such that the geometric intersection number
 between $b'(c_i)$ and $b(c_j)$ is equal to that
 between $b'(c_i)$ and $c_j$ for any $i, j=0, \ldots, n$
 and any $b'$ in the affine braid group.
 Then $b$ is the identity.

\item \label{it:fuk_faithful}
% (\ref{it:floer_gin}) and (\ref{it:gin_faithful})
 (vi) and (vii)
 suffice to prove the faithfulness
 of the action of the affine braid group
 on $D^b \Fuk W$, and hence
 on $D^b \module_0 \palg \cong \scD_k$.

\end{enumerate}

%%%%%%%%%%%%%%%%%%%%%%%%%%%%%%%%%%%%%%%%%%%%%%%%%%%%%%%%%%%
\subsection{$A_\infty$-categories and twisted complexes}

Here we recall the rudiments of $A_\infty$-categories.
For a $\bZ$-graded $k$-vector space
$N = \oplus_{j \in \bZ} N^j$ and an integer $i$,
$N[i]$ denotes its $i$-shift to the left;
$(N[i])^j  = N^{i+j}$.

\begin{definition} \label{def:a_infinity}
An $A_\infty$-category $\scA$ consists of
\begin{itemize}
 \item the set $\Ob(\scA)$ of objects,
 \item for $c_1,\; c_2 \in \Ob(\scA)$,
       a $\bZ$-graded $k$-vector space $\hom_\scA(c_1, c_2)$
       called the space of morphisms, and
 \item operations
$$\
 \m_l : \hom_\scA (c_{l-1},c_l) [1] \otimes \dots
          \otimes \hom_\scA (c_0,c_1) [1]
 \longrightarrow \hom_\scA (c_0,c_l) [1]
$$
of degree $+1$ for $l=1,2,\ldots$,
satisfying
the {\em $A_\infty$-relations}
\begin{eqnarray}
 \sum_{i=0}^{l-1} \sum_{j=i+1}^l
%  (-1)^{\deg a_1 + \cdots + \deg a_i}
  \m_{l+i-j+1}(a_l \otimes \cdots \otimes a_{j+1}
   \otimes 
    \m_{j-i}(a_j \otimes \cdots \otimes a_{i+1})
     \nonumber \\
   \otimes
    a_i \otimes \cdots \otimes a_1 ) = 0,
  \label{eq:A_infty}
\end{eqnarray}
for any positive integer $l$,
any sequence $c_0, \dots, c_l$ of objects of $\scA$,
and any sequence of morphisms
$a_i \in \hom_{\scA}(c_{i-1}, c_i)$
for $i = 1, \dots, l$.
\end{itemize}
\end{definition}

Since the $A_\infty$-relation
(\ref{eq:A_infty})
for $l=1$ ensures that
$\m_1$ squares to zero,
we can define
the {\em cohomology category}
$H^0(\scA)$
of an $A_\infty$-category $\scA$
by
$\Ob(H^0(\scA)) = \Ob(\scA)$
and
$
 \hom_{H^0(\scA)}(c_0,c_1)
  = H^0(\hom_{\scA}(c_0,c_1),\m_1).
$

To define the derived category of an $A_\infty$-category,
we need the concept of {\em twisted complexes}.
It is originally due to
Bondal and Kapranov
\cite{Bondal-Kapranov_ETC}
in the case of differential graded categories
(i.e., when $\m_k=0$ for $k \geq 3$),
and generalized to $A_\infty$-categories
by Kontsevich \cite{Kontsevich_HAMS}.

\begin{definition} \label{def:twisted-complex}
Let $\scA$ be an $A_\infty$-category.
\begin{enumerate}
 \item
The {\em additive enlargement} $\Sigma \scA$
is obtained from $\scA$
by formally adding direct sums and shifts;
an object $c$ of $\Sigma \scA$ is a formal sum
$$
 c = \bigoplus_{i \in I} c_i[l_i]
$$
where $I$ is a finite index set,
$c_i \in \Ob(\scA)$,
and
$l_i \in \bZ$.
The space of morphisms is given by
$$
 \hom_{\Sigma \scA}
  \left( \bigoplus_{i \in I} c_i[l_i],
   \bigoplus_{j \in J} d_j[m_j] \right)
  = \bigoplus_{i,j} \hom_{\scA}(c_i,d_j)[m_j-l_i],
$$
with the obvious $A_\infty$-operations
inherited from $\scA$.
%with an additional sign;
%for
%$
%a_1 \in \hom_{\scA} (c_0[l_0], c_1[l_1]),
%\dots,
%a_k \in \hom_{\scA} (c_{k-1}[l_{k-1}], c_{k}[l_k]),
%$
%$$
% \m_k(a_k,\dots,a_1)=(-1)^{l_0}(a_k,\dots,a_1).
%$$
 \item
A {\em twisted complex} over $\scA$ is a set
$(\{ c_i \}_{i \in I}, \{ x_{ij} \}_{i,j \in I})$,
where
$I \subset \bZ$ is a finite index set,
$c_i$ is an object of $\Sigma \scA$,
and
$x_{ij}$
is an element of
$\hom_{\Sigma \scA}^{i-j+1}(c_i,c_j)$
satisfying
\begin{equation*} \label{eq:twisted_complex}
 \sum_{n=1}^\infty \sum_{l=l_0<l_1<\dots<l_{n}=m}
  \m_n(x_{l_{n-1}l_n} \otimes \dots \otimes x_{l_0l_1}) = 0
\end{equation*}
for any $l, m \in \bZ$.
We assume $x_{ij} = 0$
for $i > j$.
\end{enumerate}
\end{definition}

Twisted complexes form an $A_\infty$-category:

\begin{lemma}
For an $A_\infty$-category $\scA$,
there is another $A_\infty$-category
$\PreTr(\scA)$
whose objects are
twisted complexes of $\scA$
such that the space of morphisms
between two twisted complexes
$c = (\{c_i\}_i, \{x_{ij}\}_{i,j})$ and
$d = (\{d_i\}_i, \{y_{ij}\}_{i,j})$ is
$$
 \bigoplus_{i,j}
  \hom_{\Sigma \scA}
   (c_i[i], d_j[j]).
$$
\end{lemma}

See, e.g., \cite{Fukaya_FHMS2}
for an explicit formula
of the $A_\infty$-operations
and the proof of the $A_\infty$-relations.

The cohomology category
of $\PreTr(\scA)$
is called
the {\em bounded derived category} of $\scA$
and will be denoted by $D^b(\scA)$.
%$\Ob(D^b(\scA)) = \Ob(\PreTr(\scA))$
%and
%$
% \Hom_{D^b(\scA)}(c,d)
%  = H^0(\hom_{\PreTr(\scA)}(c,d), \m_1)
%$
%for $c, d \in \Ob(D^b(\scA))$.
It has a natural structure
of a triangulated category
by a straightforward adaptation
of
\cite[proposition 2]{Bondal-Kapranov_ETC}
to the $A_\infty$ situation.

%%%%%%%%%%%%%%%%%%%%%%%%%%%%%%%%%%%%%%%%%%%%%%%%%%%%%%%%%%%
\subsection{Fukaya category}

For a pair of
a symplectic manifold $M$
and a family
$\{ L_i \}_i$
of at most countably many
Lagrangian submanifolds,
its Fukaya category
is the $A_\infty$-category
whose set of objects
is $\{ L_i \}_i$
and whose spaces of morphisms
are the Lagrangian intersection
Floer complexes.
The Floer cohomology of two exact Lagrangian submanifolds
is defined by Floer
\cite{Floer_MTLI}.
The $A_\infty$-structure
is introduced by Fukaya
\cite{Fukaya_MHACFH}
and
used by Kontsevich \cite{Kontsevich_HAMS}
to formulate his
{\em homological mirror symmetry} conjecture.
The task of defining it
in full generality
is undertaken by
Fukaya, Oh, Ohta, and Ono
\cite{Fukaya-Oh-Ohta-Ono},
although we deal only with exact Lagrangian submanifolds
in exact symplectic manifolds
in this paper.

\begin{definition}
A symplectic manifold $M$
with its symplectic form $\omega$
is exact if
there is a one-form $\theta$ on $M$
such that $\omega = d \theta$.
A Lagrangian submanifold $L$ of $M$
is exact
if there is a function $\phi$ on $L$
such that $\theta|_L = d \phi$.
\end{definition}

For two exact Lagrangian submanifolds
$L$ and $L'$
intersecting transversally,
their Floer complex
$\hom(L, L')$
is defined as the $k$-vector space
spanned by their intersection points:
$$
 \hom(L, L') = \bigoplus_{p \in L \cap L'} k[p].
$$
When the intersection is not transversal
but clean,
we can find an exact symplectomorphism
$\phi: M \to M$
so that
$L$ and $\phi(L')$ intersect transversally.
The quasi-isomorphism class of
the resulting $\hom(L, \phi(L'))$
does not depend on the choice of $\phi$
due to the basic isotopy invariance property
of the Lagrangian intersection Floer theory.

In general,
the Floer complexes do not admit $\bZ$-gradings.
To equip them with $\bZ$-gradings,
we need the following concept
of graded Lagrangian submanifolds,
introduced by Kontsevich \cite{Kontsevich_HAMS}.
%For more details, see, e.g., \cite{Seidel_GLS}.

For a symplectic manifold $M$
of dimension $2 d$
with a compatible almost complex structure $J$
such that
$2 c_1(M, J) \in H^2(M, \bZ)$ is zero,
a {\em grading} of $M$
is a choice of
(a homotopy class of)
a nowhere-vanishing section
$\Omega^{\otimes 2}$
of $(\wedge^d (T^* M, J))^{\otimes 2}$,
where $(T^* M, J)$ is the holomorphic part
of the cotangent bundle.
%Here we put $2d=\dim M$.
For a Lagrangian submanifold $L$
of $M$
with a grading $\Omega^{\otimes 2}$,
define a function
$s_L : L \to S^1 = \bCx / \bR_{>0}$ by
sending $x \in L$ to
$$
 s_L(x)
  = [\Omega^{\otimes 2}
      \left(
        (e_1 \wedge \dots \wedge e_d)^{\otimes 2}
      \right)
    ],
$$
where
$\{e_i\}_{i=1}^d$
is a basis of $T_x L$.
A {\em grading} of $L$
is a choice of a lift
$
 \stilde_L : L \to \bR
$
of $s_L$
to the universal cover $\bR$
of $S^1$.
The pair
$\Ltilde = (L, \stilde_L)$
of a Lagrangian submanifold $L$
and its grading $\stilde_L$
is called a {\em graded Lagrangian submanifold}.
For a graded Lagrangian submanifold
$\Ltilde = (L, \stilde_L)$
and an integer $l \in \bZ$,
$\Ltilde[l]$ will denote
the graded Lagrangian submanifold
$(L, \stilde_L + l)$.
This defines a free action of $\bZ$
on the set of gradings on $L$.

For an intersection point
$p \in L \cap L'$
of two graded Lagrangian submanifolds
$\Ltilde$ and $\Ltilde'$,
one can define its {\em Maslov index}
$\mu(\Ltilde, \Ltilde'; p)$.
The importance of Maslov indices
lies in the fact that
they appear in the index theorem
for Cauchy-Riemann operators
with Lagrangian boundary conditions,
and hence in the virtual dimension
of the moduli space of pseudoholomorphic maps.
See, e.g., \cite{Seidel_GLS}
for the definition
and basic properties
of Maslov indices.
For two transversal graded Lagrangian submanifolds
$\Ltilde$ and $\Ltilde'$,
we can introduce a $\bZ$-grading
on the Floer complex
$\hom(L, L')$
so that the basis $[p] \in \hom(L, L')$
for $p \in L \cap L'$
is homogeneous of degree
$\mu(\Ltilde, \Ltilde'; p)$.
We write $\hom(\Ltilde, \Ltilde')$
for the Floer complex
$\hom(L, L')$
equipped with this grading.

The $A_\infty$-operations are defined as follows:
Let
$\Ltilde_0,\ldots,\Ltilde_l$
be a sequence of
mutually transversal
graded Lagrangian submanifolds 
and $p_i \in L_i \cap L_{i+1}$,
$i = 0, \ldots, l$
be a collection of their intersection points.
Here, $L_{l+1} = L_0$ by notation.
Then the coefficient of
$p_l$ in $\m_l(p_0, \dots, p_{l-1})$
is given by the virtual number of
pseudoholomorphic maps
$\varphi: (D^2,\vec{z}) \rightarrow M$
from an $(l+1)$-pointed disk
$(D^2, \vec{z})$
such that
$\varphi(\partial_i D^2) \subset L_i$
and
$\varphi(z_i) = p_i$.
Here,
an $(l+1)$-pointed disk
is a pair
of a two-dimensional disk $D^2$
with the standard complex structure
and $l+1$ points
$\vec{z} = (z_0, \dots, z_l)$
on its boundary
respecting the cyclic order,
and
$\partial_i D^2$
denotes the interval on the boundary of $D^2$
between $z_{i-1}$ and $z_i$
for $i = 0, \dots, l$.

%%%%%%%%%%%%%%%%%%%%%%%%%%%%%%%%%%%%%%%%%%%%%%%%%%%%%%%%%%%
\subsection{Homological mirror symmetry for $A_n$-singularities}

The mirrors of $A_n$-singularities were introduced
by Hashiba and Naka \cite{Hashiba-Naka},
following the construction of
Hori, Iqbal, and Vafa
\cite{Hori-Iqbal-Vafa, Hori-Vafa}.
Our treatment here follows
Khovanov and Seidel \cite{Khovanov-Seidel} closely,
which amounts to the homological mirror symmetry
for $\scC$.
Introduce the affine manifold
$$
 W = \{ (x, y, z) \in \bC^2 \times \bCx \mid x y z = z^{n+1} - 1 \}
$$
equipped with the exact symplectic structure
$$
 \omega
  = - \frac{\sqrt{-1}}{2} \left. \left(
        d x \wedge d \overline{x} + d y \wedge d \overline{y}
        + \frac{d z \wedge d \overline{z}}{|z|^2}
       \right) \right|_W.
$$
The fiber of the projection
$$
\begin{array}{rccc}
 \pi : & W & \to & \bCx \\
 & \rotatebox{90}{$\in$} & & \rotatebox{90}{$\in$} \\
 & (x, y, z) & \mapsto & z
\end{array}
$$
is a conic in $\bC^2$,
which degenerates
to the cone
$x y = 0$
over the discriminant set
$\Delta = \{ \zeta^i \}_{i=0}^{n}$.
Here,
$\zeta = \exp[2 \pi \sqrt{-1}/(n+1)]$
is the primitive $(n+1)$th roots of unity.

A smooth map
$c : [0,1] \to \bCx$
will be called an {\em admissible curve}
if
$c(0), c(1) \in \Delta$,
$c(t) \not \in \Delta$ for $0 < t < 1$,
and
$c$ is injective.
With an admissible curve $c$,
we can associate a Lagrangian submanifold $L_c$ of $W$
by arranging the vanishing cycles of $\pi$
along $c$:
$$
 L_c = \bigcup_{0 \leq t \leq 1}
      \{(x, y, c(t)) \in W \mid |x| = |y| \}.
$$
Since $L_c$ is diffeomorphic to a sphere,
it is an exact Lagrangian submanifold.

Two admissible curves
$c$ and $c'$
are called isotopic
if there is a continuous path
$
 \phi : [0, 1] \to \Diff_0(\bCx; \Delta)
$
such that
$\phi(0) = \id$
and
$\phi(1)(c) = c'$.
We write this as
$c \simeq c'$.
In this case,
$L_{\phi(t)}$ gives a smooth family of
Lagrangian submanifolds of $W$
connecting $L_c$ and $L_{c'}$.
Since $L_c$ has the vanishing first cohomology,
it follows that
$L_c$ and $L_{c'}$ are related
by an exact symplectic isotopy.

%For a symplectic manifold $M$
%with $2 c_1(M, J) = 0$,
%there are $H^1(M, \bZ)$ ways
%of grading.
Equip $W$ with the grading
given by the second tensor power
$\Omega^{\otimes 2}$
of the holomorphic volume form
$$
 \Omega
  = \Res \frac{d x \wedge d y \wedge d z}{x y z - z^{n+1} + 1}
  = \frac{d y \wedge d z}{y z}.
$$
For an admissible curve $c$,
the Lagrangian submanifold $L_c$ of $W$
has an $S^1$-action
$$
\begin{array}{rcccc}
 S^1 \ni e^{\sqrt{-1} \theta} & : & L_c & \to & L_c \\
 & &
 \rotatebox{90}{$\in$} & &
 \rotatebox{90}{$\in$} \\
 & & (x, y, z) & \mapsto
 & (e^{\sqrt{-1} \theta} x, e^{- \sqrt{-1} \theta} y, z).
\end{array}
$$
The orbit space of this $S^1$-action
can be identified with the image
$
 c([0,1]) = \pi(L_c) \subset \bCx,
$
and for $t \in (0,1)$,
the tangent space
at a point
$(x, y, c(t)) \in L_c$
in $\pi^{-1}(c(t))$
is spanned by
%the real parts of
$
 \frac{\partial}{\partial t} = c'(t) \frac{\partial}{\partial z}
  +  \star \, \frac{\partial}{\partial y}
$
and
$
 \frac{\partial}{\partial \theta} = \sqrt{-1} y \frac{\partial}{\partial y}.
$
%Note that
%$$
% dz \Re( a \frac{\partial}{\partial z})
%  = dz \Re( (a_1 + a_2 \sqrt{-1})
%    ( \frac{\partial}{\partial z_1}
%      - \sqrt{-1} \frac{\partial}{\partial z_2} ) )
%  = (d z_1 + \sqrt{-1} d z_2)
%     (c_1 \frac{\partial}{\partial z_1}
%      + c_2 \frac{\partial}{\partial z_2})
%  = c_1 + \sqrt{-1} c_2.
%$$
Here,
we have chosen
$\frac{\partial}{\partial y}$
and
$\frac{\partial}{\partial z}$
as a basis of the tangent space of $W$,
and
$\star$ is a term which is irrelevant
for the following calculation of
$s_{L_c}$.
Then
$$
 s_{L_c}(x, y, c(t))
  = \left[ \frac{d y \wedge d z}{y c(t)}
     \left( c'(t) \frac{\partial}{\partial z}
       \wedge  \sqrt{-1} y \frac{\partial}{\partial y}
     \right)
    \right]^2
  = [- c'(t)^2 {c(t)}^{-2}].
$$
Hence
the grading of $L_c$
admits the following description:
A curve $c: [0, 1] \to \bCx$
defines a map
$s_c : [0,1] \to S^1 \cong \bCx / \bR_{>0}$
by
$
 s_c(t) = [-c'(t)^2 c(t)^{-2}],
$
and a {\em grading} of $c$
is a lift
$
 \stilde_c : [0,1] \to \bR
$
of $s_c$ to the universal cover
$\bR$ of $S^1$.
Since
$$
 s_{L_c}(x, y, c(t)) = s_c(t)
$$
for $(x, y, c(t)) \in L_c$,
a grading of $c$
is in one-to-one correspondence
with a grading of $L_c$.
%since $s_c(t) = s_{L_c}(x, y, c(t)) \in S^1$
%for any $t \in (0,1)$ and $(x, y, c(t)) \in L_c$.
For a graded curve $\ctilde$,
$\Ltilde_c$ will denote the corresponding
graded Lagrangian submanifold of $W$.
Any admissible curve admits a grading,
and there is a free and transitive action of $\bZ$
on the set of gradings on $c$.

For a transversal intersection point
$p$
of two graded curves
$\ctilde_0 = (c_0, \stilde_{c_0})$
and
$\ctilde_1 = (c_1, \stilde_{c_1})$,
define an integer
$\mu(\ctilde_0, \ctilde_1; p)$
as follows:
Fix a small circle $l$ around the point $p$
and take an intersection point
$\alpha_0 = c_0(t_0)$
of $l$ and $c_0$.
Let
$\alpha_1 = c_1(t_1)$
be the first intersection point
of $l$ and $c_1$
as one goes
from $\alpha_0$
along $l$ clockwise.
Take the clockwise arc
$\alpha : [0,1] \to l$
from $\alpha_0$ to $\alpha_1$
and a smooth map
$\pi : [0,1] \to S^1 \cong \bCx / \bR_{>0}$
such that
$\pi(0) = [-c_0'(t_0)^2 c_0(t_0)^{-2}]$,
$\pi(1) = [- c_1'(t_1)^2 c_1(t_1)^{-2}]$,
and
$\pi(t) \ne [- \alpha'(t)^2 \alpha(t)^{-2}]$
for all $t$.
Then there is a unique lift
$
 \widetilde{\pi}: [0,1] \to \bR
$
of $\pi$ such that
$
 \widetilde{\pi}(0) = \stilde_{c_0}(t_0),
$
and $\mu(\ctilde_0, \ctilde_1; p)$ is
defined by
$$
 \mu(\ctilde_0, \ctilde_1; p)
  = \stilde_{c_1}(t_1) - \widetilde{\pi}(1).
$$

Recall that two curves
$c_0$ and $c_1$ are said to intersect {\em minimally}
if they intersect transversally
and satisfy the following condition:
Take any two points $z_- \ne z_+$ in $c_0 \cap c_1$
which do not both lie in $\Delta$,
and two arcs $\alpha_0 \subset c_0$,
$\alpha_1 \subset c_1$ with endpoints $z_-$, $z_+$,
such that $\alpha_0 \cap \alpha_1 = \{ z_-, z_+\}$.
Let $K$ be the connected component of
$\bCx \setminus (c_0 \cup c_1)$
which is bounded by $\alpha_0 \cup \alpha_1$.
Then if $K$ is topologically an open disk,
it must contain at least one point of $\Delta$.

For two graded curves
$\ctilde_0$ and $\ctilde_1$
intersecting minimally,
their {\em graded intersection number}
is defined by
$$
 \Igr(\ctilde_0, \ctilde_1)
  = (1+q) \sum_{p \in (c_0 \cap c_1) \setminus \Delta}
     q^{\mu(c_0, c_1; p)}
    + \sum_{p \in c_0 \cap c_1 \cap \Delta}
     q^{\mu(c_0, c_1; p)}.
$$
Here,
the factor of $1 + q$
in front of the sum over
$(c_0 \cap c_1) \setminus \Delta$
comes from the Poincar\'{e} polynomial
$
 \sum_r q^r \dim_k H^r(S^1;k) = 1 + q
$
of $S^1$
and the fact that
$L_{c_0}$ and $L_{c_1}$ intersect cleanly
along $S^1$
over $p \in (c_0 \cap c_1) \setminus \Delta$.

The following result is due to Khovanov and Seidel:

\begin{theorem}[{\cite[lemma 6.19]{Khovanov-Seidel}}]
\label{theorem:Khovanov-Seidel}
Let $\Ltilde_0$ and $\Ltilde_1$ be the exact graded Lagrangian
two-spheres in $W$
coming from two admissible graded curves
$\ctilde_0$ and $\ctilde_1$ on the $z$-plane.
Assume that $c_0$ and $c_1$ intersect
transversally and minimally.
Then the Poincar\'{e} polynomial of the Floer cohomology groups
of $\Ltilde_0$ and $\Ltilde_1$
is given by the graded intersection number
of $\ctilde_0$ and $\ctilde_1$:
$$
 \sum_{r \in \bZ} q^r \dim_k
  H^r(\hom(\Ltilde_0, \Ltilde_1), \m_1)
   = \Igr(\ctilde_0, \ctilde_1).
$$
\end{theorem}

Since
$
 \Igr(\ctilde_0, \ctilde_1)
$
at $q=1$
does not depend on the choice of
the grading,
we write
$$
 I(c_0, c_1) = \frac{1}{2} \Igr(\ctilde_0, \ctilde_1) |_{q = 1}
$$
and call it the {\em geometric intersection number}.

\begin{figure}[htbp]
\centering
\psfrag{0}{$0$}
\psfrag{z_0}{$1$}
\psfrag{z_1}{$\zeta$}
\psfrag{z_2}{$\zeta^2$}
\psfrag{z_n-1}{$\zeta^{n-1}$}
\psfrag{z_n}{$\zeta^{n}$}
\psfrag{L_0}{$c_0$}
\psfrag{L_1}{$c_1$}
\psfrag{L_n-1}{$c_{n-1}$}
\psfrag{L_n}{$c_n$}
\includegraphics[width=0.4 \textwidth]{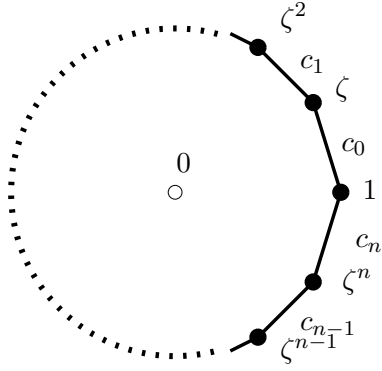}
\caption{The curves $c_i$'s}
\label{fg:vanishing_cycles}
\end{figure}

For $i = 0, \ldots, n$,
let $c_i : [0,1] \to \bCx$
be the straight line segment
from $\zeta^i$ to $\zeta^{i+1}$
as in Figure \ref{fg:vanishing_cycles}.
Equip $c_i$ with the grading
such that the value of
$\stilde_{c_i}$
at the midpoint $(\zeta_i + \zeta_{i+1}) / 2$
is the same for all $i=0, \dots, n$,
and
let $\Ltilde_i$ be
the corresponding graded Lagrangian submanifold of $W$.
Let further
$\Fuk W$
be the Fukaya category of $W$
whose set of objects
is $\{ \Ltilde_i \}_{i=0}^n$.

The following theorem is the homological mirror symmetry
for $A_n$-singularities.
Note that we do not need the assumption on the characteristic
of $k$ in the proof.

\begin{theorem}\label{theorem:HMS}
There is an equivalence
$$
  \scD_k \cong D^b \Fuk W
$$
of triangulated categories.
\end{theorem}

\begin{proof}
We can see that
there are no non-constant holomorphic maps
from a disk to $W$
with $\{ L_i \}_{i=0}^n$
as Lagrangian boundary conditions:
If $\phi: D^2 \to W$
is a holomorphic map
from a disk such that
$\phi(\partial D^2) \subset \bigcup_{i=0}^n L_i$,
then $\pi \circ \phi$ is a holomorphic map
from $D^2$ to $\bCx$
with the boundary condition
$\phi(\partial D^2) \subset \bigcup_{i=0}^n c_i$.
Since non-constant holomorphic maps are open
and there are no continuous open maps
from $D^2$ to $\bCx$
with the above boundary condition,
$\pi \circ \phi$ must be constant.
Hence $\phi$ is a holomorphic map
from $D^2$
into the fiber $\pi^{-1}(z)$
for some $z \in \bigcup_{i=0}^n c_i \subset \bCx$.
It is obvious that
for any $z \in \bCx$,
there are no non-constant holomorphic maps
$\phi : D^2 \to \pi^{-1}(z)$
such that
$$
 \phi(\partial D^2) \subset \{ (x, y, z) \in \pi^{-1}(z)
  \mid |x| = |y| \}.
$$
Therefore
the structure of the Fukaya category $\Fuk W$
is easy to describe:
The spaces of morphisms are given by
$$
 \hom_{\Fuk W} (\Ltilde_i, \Ltilde_j) =
  \begin{cases}
   k \cdot e_i \oplus k \cdot f_i & \text{if } i=j, \\
   k \cdot x_i & \text{if } j = i + 1, \\
   k \cdot y_i & \text{if } j = i - 1, \\
   0 & \text{otherwise},
  \end{cases}
$$
where $\deg(e_i) = 0$,
$\deg(x_i) = \deg(y_i) = 1$,
and $\deg(f_i) = 2$.
The $A_\infty$-operations are trivial
except for $\m_2$,
and $e_i$ is the identity morphism of $\Ltilde_i$
for $i = 0, \dots, n$.
Nontrivial $\m_2$ are given by
$$
\begin{array}{ccccc}
 \m_2(y_{i+1}, x_i) & = & \m_2(x_{i-1}, y_i) & = & f_i,
 \qquad i = 0, \dots, n. \\
\end{array}
$$
Hence
the total morphism algebra
$
 \bigoplus_{i,j=0}^n \hom_{\Fuk W}(\Ltilde_i, \Ltilde_j)
$
is not only an $A_\infty$-algebra
but a differential graded algebra
(i.e., $\m_n = 0$ for $n \geq 3$)
with the trivial differential
(i.e., $\m_1 = 0$).
Note that for any $i, j = 0, \dots, n$,
we have
$$
 \hom_{\Fuk W}(\Ltilde_i, \Ltilde_j)
  = \Ext_{\palg}^*(S_i, S_j),
$$
where $\palg$ is the path algebra with relations
appearing in \S \ref{section:McKay}
and $S_i$ is the simple $\palg$-module
corresponding to the idempotent
$(i) \in \palg$.
Since the abelian category
$\module \palg$
of finitely generated right $\palg$-modules
has enough injectives,
its derived category
$D^b \module \palg$
has an enhancement
in the sense of Bondal and Kapranov
\cite{Bondal-Kapranov_ETC}.
Since
$\palg$ is Koszul
over the semisimple ring
$k^{\oplus (n+1)}$,
the total endomorphisms $DG$-algebra
of $\bigoplus_{i=0}^n S_i$
in the enhancement of $D^b \module \palg$
is formal.
Then a theorem of Bondal and Kapranov
\cite[theorem 1]{Bondal-Kapranov_ETC}
shows that
$D^b \Fuk W$ is equivalent
as a triangulated category
to the smallest triangulated subcategory
of $D^b \module \palg$
containing $\{ S_i \}_{i=0}^n$.
This subcategory is equivalent
to the bounded derived category
$D^b \module_0 \palg$
of finitely-generated nilpotent $\palg$-modules.
Since $D^b \module_0 \palg$
is triangle-equivalent to $\scD_k$
by the McKay correspondence,
we obtain an equivalence $\scD_k \cong D^b \Fuk W$
as desired.
\end{proof}

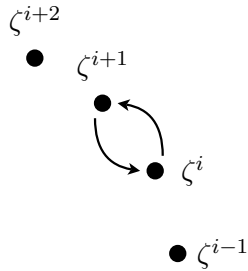
\begin{figure}[htbp]
\centering
\ifx\JPicScale\undefined\def\JPicScale{1}\fi
\psset{unit=\JPicScale mm}
\psset{linewidth=0.3,dotsep=1,hatchwidth=0.3,hatchsep=1.5,shadowsize=1,dimen=middle}
\psset{dotsize=0.7 2.5,dotscale=1 1,fillcolor=black}
\psset{arrowsize=1 2,arrowlength=1,arrowinset=0.25,tbarsize=0.7 5,bracketlength=0.15,rbracketlength=0.15}
\begin{pspicture}(0,0)(36,33)
\rput{0}(11,28){\psellipse[fillstyle=solid](0,0)(1,-1)}
\rput{0}(20,22){\psellipse[fillstyle=solid](0,0)(1,-1)}
\rput{0}(27,13){\psellipse[fillstyle=solid](0,0)(1,-1)}
\rput{0}(30,2){\psellipse[fillstyle=solid](0,0)(1,-1)}
\psbezier{<-}(22,22)(26.2,21.3)(28,19.2)(28,15)
\psbezier{->}(19,20)(19,15.8)(20.8,13.7)(25,13)
\rput(32,13){$\zeta^i$}
\rput(36,2){$\zeta^{i-1}$}
\rput(20,27){$\zeta^{i+1}$}
\rput(11,33){$\zeta^{i+2}$}
\end{pspicture}
\caption{The Dehn twist $\tau_i$}
\label{fg:dehn_twist}
\end{figure}

Now we discuss the symplectic Dehn twist.
Let
$$
 T = \{ (a_0, \dots, a_n) \in \bC^{n+1} \mid a_0 \cdots a_n = 1 \}
$$
be an algebraic torus of dimension $n$
and
$$
 D = \{ (a_0, \dots, a_n) \in T \mid
  a_i = a_j \ \text{for some} \ i \neq j \}
$$
be its big diagonal.
The symmetric group $\frakS_{n+1}$
of rank $n+1$ acts
on $T$ by permutations,
and let $S$ be the quotient
$(T \setminus D) / \frakS_{n+1}$.
Consider the family
$$
 \scW = \{
  (x, y, z, a) \in \bC^2 \times \bCx \times S
   \mid
  x y z = (z - a_0) \cdots (z - a_n)
  \} \to S
$$
of exact symplectic manifolds over $S$,
equipped with the relative grading
given by the relative holomorphic volume form
$\Omega = d y / y \wedge d z / z$.
Since any fibration of exact symplectic manifolds
is locally trivial by Moser \cite{Moser_OVEM},
a loop in the parameter space $S$
induces a graded symplectomorphism
of the fiber as its monodromy.
For $i = 0, \dots, n$,
let $\tau_i$ be the monodromy
along the loop
around $a_i = a_{i+1}$
as in Figure \ref{fg:dehn_twist},
starting and ending at
$
 [(a_0, a_1, \dots, a_n)]
  = [(1, \zeta, \dots, \zeta^n)]
 \in S.
$
Then $\tau_i$ acts
both on graded curves
and on graded Lagrangian submanifolds,
and the construction of graded Lagrangian submanifolds
from graded curves
commutes with this action;
for a graded curve $\ctilde$,
we have
$$
 \tau_i(\Ltilde_{\ctilde})
  = \Ltilde_{\tau_i(\ctilde)}.
$$
Let $\scFtilde$ be the Fukaya category
whose set of objects
is the images of
$L_i$,
$i=0, \dots, n$
by the compositions of
$\tau_j$
for
$j=0, \dots, n$.
The following theorem is due to Seidel:
\begin{theorem}[{\cite[proposition 9.1]{Seidel_K3}}]\label{th:Seidel}
The graded Lagrangian submanifold
$\tau_i(\Ltilde)$
is isomorphic to the twisted complex
$T_{\Ltilde_i}(\Ltilde)$
in $D^b \scFtilde$.
\end{theorem}
The above theorem also shows that
$D^b \scFtilde$ is equivalent to $D^b \Fuk W$
as a triangulated category.

%%%%%%%%%%%%%%%%%%%%%%%%%%%%%%%%%%%%%%%%%%%%%%%%%%%%%%%%%%%
\subsection{The proof of the faithfulness in characteristic two}

Here we prove the faithfulness
of the  affine braid group action
on $D^b \scFtilde$
given by
$$
\begin{array}{rcccc}
 \rho & : & \ABG & \to & \Auteq D^b \scFtilde \\
 & &
 \rotatebox{90}{$\in$} & &
 \rotatebox{90}{$\in$} \\
 & & \sigma_i & \mapsto
 & \tau_i.
\end{array}
$$
Since the action of $\ABG$
on graded curves
and graded Lagrangian submanifolds
commute,
%$$
% \Phi_b(\Ltilde_{\ctilde}) \cong \Ltilde_{b(\ctilde)},
%$$
we can work with graded curves on $\bCx$
instead of graded Lagrangian submanifolds of $W$.

The following lemma is an affine version
of \cite[lemma 3.6]{Khovanov-Seidel}:

\begin{lemma} \label{lem:faithfulness}
If $b \in \ABG$ satisfies
\begin{equation} \label{eq:lemma_faithfulness}
 I(b(c_i), b'(c_j)) = I(c_i, b'(c_j))
\end{equation}
for any $i, j = 0, \ldots, n$
and any $b' \in \ABG$,
then $b$ is the identity.
\end{lemma}

\begin{proof}
Think of $\ABG$
as a subgroup of $\pi_0(\Diff_0(\bCx; \Delta))$
and take an element $b$ as above.
We prove that
$b(c_i) \simeq c_i$ for $i=0, \ldots, n$,
which implies $b = \id$.
It suffices to prove $b(c_0) \simeq c_0$
by symmetry.
Put $c = b(c_0)$.
By substituting $b' = \id$,
$i = 0$, and $j = 1, \dots, n$
into (\ref{eq:lemma_faithfulness}),
we obtain
$$
 I(c, c_1) = I(c, c_n) = \frac{1}{2}
$$
and
$$
 I(c, c_2) = \cdots = I(c, c_{n-1}) = 0.
$$
This leaves two possibilities
$c'$ and $c''$
in Figure \ref{fg:possible_c}
for the homotopy class of $c$
other than $c_0$.
\begin{figure}[htbp]
\centering
\psfrag{0}{$0$}
\psfrag{c'}{$c'$}
\psfrag{c''}{$c''$}
\psfrag{c10}{$\tau_1 c_0$}
\psfrag{c-10}{$\tau_1^{-1} c_0$}
\psfrag{z_0}{$1$}
\psfrag{z_1}{$\zeta$}
\includegraphics[width = .5 \textwidth]{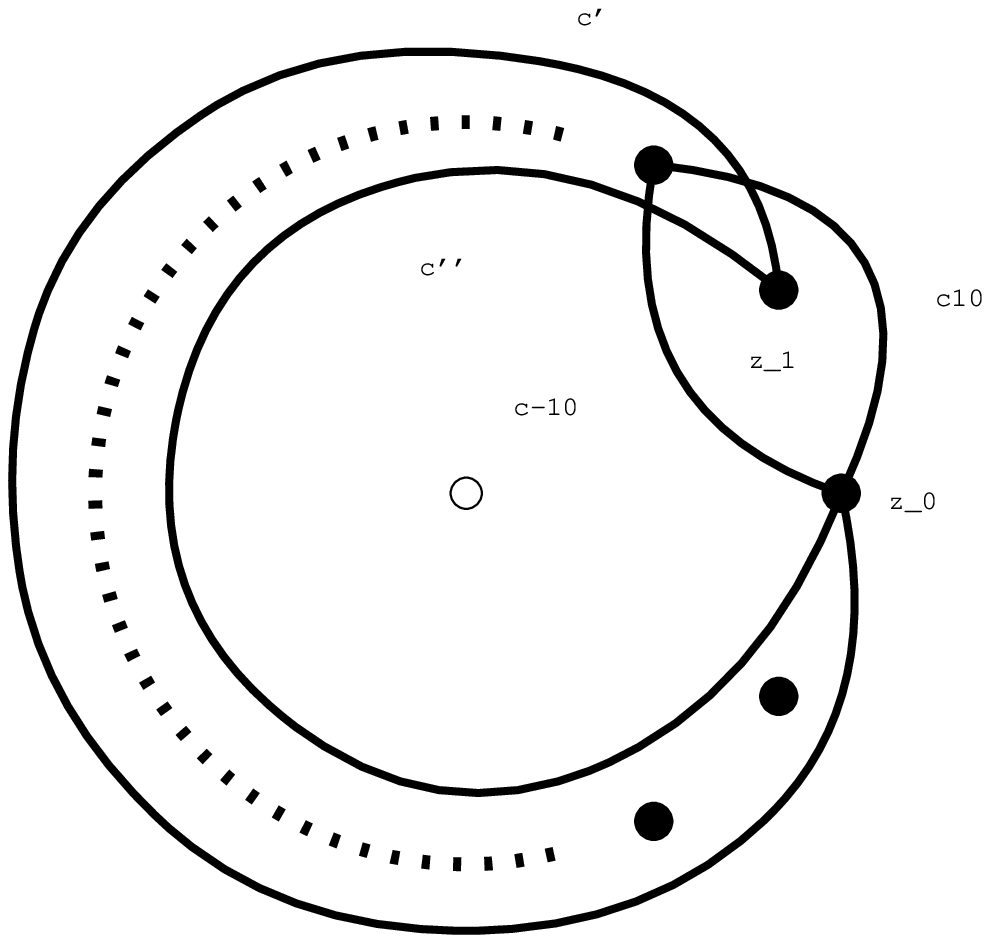}
\caption{$c'$ and $c''$}
\label{fg:possible_c}
\end{figure}
However,
we have
$$
 I(c', \tau_1 c_0)
    = 1
 \neq \frac{1}{2}
    = I(c_0, \tau_1 c_0)
$$
and
$$
 I(c'', \tau_1^{-1} c_0)
    = 1
 \neq \frac{1}{2}
    = I(c_0, \tau_1^{-1} c_0),
$$
which shows that
neither
$c'$ nor $c''$ satisfies
(\ref{eq:lemma_faithfulness})
for $b' = \sigma_1$ and $b' = \sigma_1^{-1}$
respectively.
Hence we have
$c \simeq c_0$.
\end{proof}

Now we can prove the faithfulness in characteristic two:

\begin{theorem} \label{th:faithfulness}
If $b \in \ABG$ satisfies
$\rho(b) = \id$ in $\Auteq D^b \scFtilde$,
then $b = \id $ in $\ABG$.
\end{theorem}

\begin{proof}
Since $\rho(b) = \id$,
$$
 \dim \Hom^*(\rho(b)(L_i), \rho(b')(L_j))
  = \dim \Hom^*(L_i, \rho(b')(L_j))
$$
for any $i, j = 0, \dots, n$ and
any $b' \in \ABG$.
Since
$$
 \dim \Hom^*(\rho(b)(L_i), \rho(b')(L_j))
  = 2 I(b(c_i), b'(c_j)),
$$
by Theorem \ref{theorem:Khovanov-Seidel},
we obtain $b = \id$
from Lemma \ref{lem:faithfulness} .
\end{proof}

%%%%%%%%%%%%%%%%%%%%%%%%%%%%%%%%%%%%%%%%%%%%%%%%%%%%%%%%%%%%%%%%%%%%%%%%
%%%%%%%%%%%%%%%%%%%%%%%%%%%%%%%%%%%%%%%%%%%%%%%%%%%%%%%%%%%%%%%%%%%%%%%%
\section{Lifting to any characteristic} \label{section:lifting}

In this section,
we lift the above result to any characteristic.
For that purpose, we consider
$$
Y_{\bZ}= \Spec \bZ[x, y, z]/(xy+z^{n+1})
$$
and the minimal resolution
$f: X_{\bZ} \to Y_{\bZ}$
with the exceptional locus $Z_{\bZ}=\cup_{i=1}^n C_{i, \bZ}$.
For any noetherian ring $R$, we denote $X_R=X_{\bZ} \otimes_{\bZ} R$,
$Y_R=Y_{\bZ} \otimes_{\bZ} R$,
$\scD_R=D^b \coh_{Z_R} X_R$, etc.

We say that an object $\scE \in \scD_R$ is spherical if
\begin{enumerate}
\item $\scE$ is quasi-isomorphic to a bounded complex of $R$-flat 
coherent sheaves
on $X_R$ whose cohomology sheaves are supported at $Z_R$,
\item
$\Ext^0_{X_R}(\scE, \scE) \cong R \cong \Ext^2_{X_R}(\scE, \scE)$,
\item
$\Ext^i_{X_R}(\scE, \scE)=0$ for $i \ne 0, 2$.
\end{enumerate}

\begin{remark}
\begin{enumerate}
\item
Since $X_{\bZ}$ is smooth over $\bZ$, we can replace ``$R$-flat" in the 
first condition
by ``locally free."
\item
When $R$ is regular, then the first condition is automatically 
satisfied.
\end{enumerate}
\end{remark}

For an object $\scE \in \scD_R$ and
any noetherian $R$-algebra $S$, we put $\scE_S:=\scE 
\overset{\bL}{\otimes}_R S$.
 
\begin{lemma}
If $\scE \in \scD_R$ is spherical,
so is $\scE_S$.
\end{lemma}
\begin{proof}
This follows from the base change theorem.
(See \cite[\S 7.7 and \S 7.8]{EGAIII}.)
\end{proof}

\begin{lemma}
Let $(R, m)$ be a noetherian local ring with residue field $k'$.
Suppose that $\scE_0 \in \scD_{k'}$ is a spherical object.
Then there is a spherical object $\scE \in \scD_R$
with $R$-flat cohomology sheaves $\scH^i(\scE)$
which satisfies
$\scE_{k'} \cong \scE_0$.
\end{lemma}
\begin{proof}
Put $\scH_0^i:=\scH^i(\scE_0)$.
Then, as in \cite[\S 4.1 and \S 4.2]{Ishii-Uehara_ADC},
$\scE_0$ determines elements $e^i(\scE_0) \in \Ext^2_{X_{k'}}(\scH_0^i, 
\scH_0^{i-1})$
and the isomorphism class of $\scE_0$ is determined by these data.
We construct $\scE$ by lifting these data to $X_R$.
$\scH_0^i$ is a direct sum
of sheaves that are line bundles on subchains of exceptional curves
and hence can be lifted to a sheaf $\scH^i \in \coh_{Z_R} X_R$ flat 
over $R$.
Since $\Ext^p_{X_{k'}}(\scH_0^i, \scH_0^j)=0$ for $p \ne 0, 2$ by 
\cite[proposition 4.5]{Ishii-Uehara_ADC},
the base change theorem implies that $\Ext^p_{X_R}(\scH^i, \scH^j)$ is 
$R$-free and
$\Ext^p_{X_{k'}}(\scH_0^i, \scH_0^j) \cong \Ext^p_{X_R}(\scH^i, \scH^j) 
\otimes_R k'$.
Especially, $e^i(\scE_0)$ can be lifted to an element
$e^i \in \Ext^2_{X_R}(\scH^i, \scH^{i-1})$.
As in \cite[proposition 4.2]{Ishii-Uehara_ADC}, it follows from $\Ext^{\ge 
3}_{X_R}(\scH^i, \scH^j)=0$ that there is an object $\scE \in \scD_R$
with $\scH^i(\scE) \cong \scH^i$ and $e^i(\scE)=e^i$, which is unique 
up to isomorphism.
Then the construction yields $\scE_{k'} \cong \scE_0$ and therefore $\scE$ 
is spherical
by the base change theorem.
\end{proof}

\begin{corollary}\label{corollary:flat}
Let $\scE \in \scD_R$ be spherical for some noetherian ring $R$.
Then the cohomology sheaves $\scH^i(\scE)$ are $R$-flat.
\end{corollary}
\begin{proof}
We may assume $R$ is local with residue field $k'$.
By the previous lemma,
there is a spherical object $\scF \in \scD_R$
with $R$-flat cohomology sheaves
satisfying $\scF_{k'} \cong \scE_{k'}$.
Since we have $\Ext^1_{X_{k'}}(\scE_{k'}, \scE_{k'}) =0$ and 
$\Ext^0_{X_{k'}}(\scE_{k'}, \scE_{k'}) \cong {k'}$,
the deformation theory in \cite{Inaba_TDMCCS} implies that $\scE$ and $\scF$ 
are isomorphic
on some \'etale covering of $\Spec R$.
This proves that $\scH^i(\scE)$ is $R$-flat.
\end{proof}

For a spherical object $\scE \in \scD_R$,
we can define the twist functor $T_{\scE}$ as a Fourier-Mukai transform
with respect to $\{ \scE^{\vee} \boxtimes \scE \to \scO_{\Delta}\}$.
Then, for any $R$-algebra $S$, the base change theorem (see 
\cite[\S 2.1]{Bonsdorff_FTHB}) implies that the twist functors $T_{\scE}$ and 
$T_{\scE_S}$
commute with the functor $\cdot\, \overset{\bL}{\otimes}_R S: \scD_R 
\to \scD_S$.
In particular,
for a spherical object $\scE$ in $\scD_{\bZ}$, we have the 
following commutative diagram:
$$
\begin{CD}
\scD_{k} @<<< \scD_{\bZ} @>>> \scD_{\bF_2} \\
@V{T_{\scE_{k}}}VV        @V{T_{\scE}}VV  @V{T_{\scE_{\bF_2}}}VV\\
\scD_{k} @<<< \scD_{\bZ} @>>> \scD_{\bF_2}
\end{CD}
$$
Recall that $k$ is the base field of any characteristic.
Thus for an element $b$ of the affine braid group, we have 
autoequivalences
$\rho_{k}(b)$, $\rho_{\bZ}(b)$, $\rho_{\bF_2}(b)$
of $\scD_{k}$, $\scD_{\bZ}$, $\scD_{\bF_2}$,
respectively, which commute with
$\cdot\, \overset{\bL}{\otimes}_{\bZ} k$ and 
$\cdot\, \overset{\bL}{\otimes}_{\bZ} \bF_2$.

\begin{proposition}\label{prop:lifting}
If $\rho_{k}(b)$ is isomorphic to the identity,
then $\rho_{\bF_2}(b)$ is also isomorphic to the identity.
\end{proposition}

\begin{proof}
%It suffices to show that
%$\rho_{\bF_2}(b)(\scO_{C_{i,\bF_2}}(d)) \cong 
%\scO_{C_{i,\bF_2}}(d)$
%for any sheaf of the form $\scO_{C_{i, \bF_2}}(d)$.
Since $\rho_{k}(b)$ is the identity, we have
$\rho_{\bZ}(b)(\scO_{C_{i, \bZ}}(d)) \otimes k \cong \scO_{C_{i, k}}(d)$
for any $d \in \bZ$.
Now $\rho_{\bZ}(b)(\scO_{C_{i, \bZ}}(d))$ has flat cohomologies by 
Corollary \ref{corollary:flat} and hence
we have $\rho_{\bZ}(b)(\scO_{C_{i, \bZ}}(d)) \cong \scO_{C_{i, \bZ}}(d)$.
This implies $\rho_{\bF_2}(b)(\scO_{C_{i,\bF_2}}(d)) \cong 
\scO_{C_{i,\bF_2}}(d)$.
Therefore $\rho_{\bF_2}(b)$ is a Fourier-Mukai functor
which sends
the structure sheaf of a closed point
to the structure sheaf of a closed point.
Since $\rho_{\bF_2}(b)$ acts as the identity functor
on objects supported outside the exceptional set,
$\rho_{\bF_2}(b)$ is isomorphic
to a functor of the form
$\cdot \otimes \scL$
for some line bundle $\scL$ on $X$.
The fact that
$\rho_{\bF_2}(b)(\scO_{C_i, \bF_2}(d)) \cong \scO_{C_i, \bF_2}(d)$
for any $d$
shows that $\scL$ must be trivial.
\end{proof}
By combining Theorems \ref{theorem:HMS}, \ref {th:Seidel},
\ref{th:faithfulness} and
Proposition \ref{prop:lifting},
we obtain the following:
\begin{corollary}
The homomorphism $\rho_k$ is injective
for any field $k$.
\end{corollary}

%%%%%%%%%%%%%%%%%%%%%%%%%%%%%%%%%%%%%%%%%%%%%%%%%%%%%%%%%%%%%%%%%%%%%%%%
%%%%%%%%%%%%%%%%%%%%%%%%%%%%%%%%%%%%%%%%%%%%%%%%%%%%%%%%%%%%%%%%%%%%%%%%
\section{Connectedness of $\Stab \scC$} \label{section:connectedness_C}

In this section,
we prove the connectedness of $\Stab \scC$.
Our strategy is the following:
\begin{enumerate}
 \item First we prove
Lemma \ref{lemma:length},
which states
that any spherical object of $\scC$ can be obtained
from $\scO_{C_i}(-1)$ for some $i=1, \dots, n$
by the action of $\Br(\scC)$.
This shows, on the mirror side,
that any spherical object can be represented
not only by a twisted complex of graded Lagrangian submanifolds
but by an honest graded Lagrangian submanifold of $W$
obtained by twisting $L_1, \dots, L_n$ along themselves.
This opens up a way to use topological arguments
on configurations of curves on a disk
to tackle the problem of the connectedness of $\Stab \scC$.
 \item Let $\Sigma$ be a connected component
of $\Stab \scC$.
By Lemma \ref{lemma:perturb},
we can find a stability condition
$\sigma = (Z, \scP) \in \Sigma$
such that $Z(E) \not \in \bR$
for any spherical object $E \in \scC$.
 \item Let $\{ E_1, \dots, E_m \}$
be the set of simple objects
of the heart $\scP((0,1])$
of the above stability condition.
Then $m=n$ and
$\{ E_1, \dots, E_n \}$ forms an $A_n$-configuration
with a suitable order.
 \item The homological mirror symmetry
allows us to use
a topological argument
on a configuration of curves on a disk
to find an autoequivalence $\Phi \in \Br(\scC)$
which brings $E_i$ to the standard generator
$\scO_{C_i}(-1)$ simultaneously
for all $i = 1, \dots, n$
in characteristic two.
Then the lifting argument in \S \ref{section:lifting}
allows us to prove it
in any characteristic.
\item The above autoequivalence $\Phi$ brings $\sigma$
to $\Phi \sigma$,
which belongs to the distinguished connected
component of $\Stab \scC$
studied by Bridgeland.
Since the action of $\Br(\scC)$
preserves this connected component
\cite[theorem 1.1]{Bridgeland_SCKS},
$\Stab \scC$ is connected.
\end{enumerate}

%%%%%%%%%%%%%%%%%%%%%%%%%%%%%%%%%%%%%%%%%%%%%%%%%%%%%%%%%%%
\subsection{Normalizing spherical objects}

The following lemma follows from
the result in \cite{Ishii-Uehara_ADC}:

\begin{lemma}\label{lemma:length}
Let $E$ be a spherical object in $\scC$.
Then there is an autoequivalence $\Phi\in \Br(\scC)$
such that $\Phi (E)\cong \scO_{C_i}(-1)$ for some $i$.
\end{lemma}

\begin{proof}
First note the fact that
$E\in\scC$ if and only if $\scH ^j(E)\in\scC$ for any $j$
by \cite[lemma 3.1]{Bridgeland_FDC}.
Combining this with
the proof of Proposition 1.6 in
\cite{Ishii-Uehara_ADC},
we can show that for a spherical object $E\in \scC$
with $l(E)>1$,
there is an autoequivalence $\Phi\in \Br(\scC)$
such that $l(E)>l(\Phi(E))$.
Here,
the length $l(E)$ of an object $E \in \scD$ is defined
in \cite{Ishii-Uehara_ADC} as
$$
l(E)=\sum_{i,p}\operatorname{length}_{\scO_{X,\eta_i}} \scH^p(E)_{\eta_i},
$$
where $\eta_i$ is the generic point of $C_i$.
We can check that the spherical object $E\in\scC$ with $l(E)=1$
has the form $\scO_{C_i}(-1)[d]$ for some $i$, $d\in\bZ$.
Hence the equality
$T_{\scO_{C_i}(-1)}(\scO_{C_i}(-1))=\scO_{C_i}(-1)[-1]$
completes the proof.
\end{proof}
 
%%%%%%%%%%%%%%%%%%%%%%%%%%%%%%%%%%%%%%%%%%%%%%%%%%%%%%%%%%%
\subsection{$A_n$-configurations from stability conditions}

Fix a connected component $\Sigma$ of $\Stab\scC$.
We can show $V(\Sigma)=\Hom(K(\scC), \bC)$
in just the same way as in
the proof of Lemma \ref{lemma:perturb}.
Hence we have a stability condition $\sigma=(Z, \scP)\in \Sigma$
such that $Z(E)\not\in\bR$ for any spherical object $E \in \scC$.
Then $\Im Z$ determines a Weyl chamber
and the corresponding simple root basis.
We prove that
for such a stability condition,
the set of all simple objects in its heart $\scP((0,1])$
forms an $A_n$-configuration
whose image in $K(\scC)$ is the simple root basis.

%%%
%%%
%%%

\begin{lemma}\label{lemma:non-degenerate}
Let $\sigma = (Z, \scP) \in \Sigma$ be a stability condition
such that $Z(E)\not\in\bR$ for any spherical object $E \in \scC$.
Then the set of all simple objects in $\scP((0,1])$
consists of $n$ mutually non-isomorphic elements. 
%and $\bigl\{E_1,\ldots,E_m \bigr\}$ be the set of all 
%simple objects in its heart $\scP((0,1])$.
%Then we have $m=n$
Moreover,
this set forms an $A_n$-configuration of spherical objects
in an appropriate order.
% i.e.
%after replacing the suffix of $E_i$'s if necessary,
%we have 
%\[
%\dim\Hom _{\scC}^*(E_i,E_j)=
%\begin{cases}
%1 & (|i-j|=1) \\
%0 & (|i-j|\ne 1).
%\end{cases}
%\]
\end{lemma}

\begin{proof}
Note that a simple object in $\scP((0,1])$ is always stable and
any stable object in $\scC$ is spherical
by Lemma \ref{lemma:stable_real_root_spherical}.
Moreover, the class of a stable object
in $\scP((0,1])$ is a positive root.
Thus $\scP((0,1])$ is a finite length category.

If two simple objects $E_1, E_2 \in \scP((0,1])$ determine the same class
in $K(\scC)$, then $\chi(E_1, E_2)=2$ holds and we have $E_1 \cong E_2$.
Hence, the set of simple objects in $\scP((0,1])$ injects into the set of
positive roots; in particular, it is finite.
Let $\bigl\{E_1,\ldots,E_m \bigr\}$ be the set of all simple objects
in $\scP((0,1])$.
Then $\scP((0,1])$ 
is the smallest extension closed subcategory of $\scC$ 
containing this set,
since it is of finite length.

If $i \ne j$,
then
we have $\Hom _{\scC}^{\le 0}(E_i,E_j) = 0$
since $E_i$ and $E_j$ are mutually 
non-isomorphic simple objects
of the heart.
Then
the Serre duality implies $\Hom _{\scC}^{\ge 2}(E_i,E_j) = 0$.
Hence
we have
\begin{equation}\label{equation:ge0}
\chi(E_i,E_j)=-\dim\Hom _{\scC}^1(E_i,E_j) \le 0.
\end{equation}
First we show that $\bigl\{ E_1, \ldots, E_m \bigr\}$ is a linear basis in $K(\scC)\otimes  
\bQ$.
Assume that we have an equality $\sum _i a_i[E_i]=0$ for some  
$a_i\in\bQ$.
Put
\[\alpha=\sum_{a_i>0}a_i[E_i]=-\sum_{a_j<0}a_j[E_j]. \]
Then (\ref{equation:ge0}) implies
\[\chi(\alpha,\alpha) =-\sum_{a_i>0,a_j<0}a_ia_j\chi(E_i,E_j)\le 0,\]
and we obtain $\alpha=0$.
Since $[E_i]$'s are positive roots, we obtain $a_i =0$ for all $i$.
Thus they form a basis of $K(\scC)\otimes{\bQ}$.
In particular we have $m=n$.

To show that
$\bigl\{ E_1, \ldots, E_m \bigr\}$
forms a simple root basis of $K(\scC)$,
let $\alpha \in K(\scC)$ be any root and write
$\alpha = \sum_{i = 1}^n a_i [E_i]$.
Put
$\alpha_1=\sum_{a_i>0}a_i[E_i]$ and $\alpha_2=\sum_{a_j<0}a_j[E_j]$.
Then we can use \eqref{equation:ge0} to see
$$
\chi(\alpha, \alpha) = \chi(\alpha_1, \alpha_1) + \chi(\alpha_2,
\alpha_2) + 2\chi(\alpha_1, \alpha_2) \ge \chi(\alpha_1, \alpha_1) +
\chi(\alpha_2, \alpha_2).
$$
Since $\alpha$ is a root, either $\alpha_1$ or $\alpha_2$ is zero.
Hence
$\bigl\{ E_1, \ldots, E_m \bigr\}$
forms a simple root basis.
%These facts ensure that $\bigl\{[E_1],\ldots,[E_m] \bigr\}$
%is actually the simple root basis.
% This is true at least in A_n case. 
Consequently
the set $\bigl\{E_1,\ldots,E_n\bigr\}$
with a suitable order
forms an $A_n$-configuration of spherical objects.
\end{proof}

%%%%%%%%%%%%%%%%%%%%%%%%%%%%%%%%%%%%%%%%%%%%%%%%%%%%%%%%%%%              
\subsection{Normalizing an $A_n$-configuration}

Here we prove the following lemma:

\begin{lemma} \label{lemma:configuration}
Let $(E_1, \dots, E_n)$
be an $A_n$-configuration of spherical objects 
in $\scC_k$,
and assume that for any $i$,
there is an element $\Phi \in \Br(\scC_k)$
such that
$E_i \cong \Phi(\scO_{C_j}(-1))$
for some $j$.
Then there is an element
$\Psi \in \Br(\scC_k)$
such that
$\Psi(E_i) \cong \scO_{C_i}(-1)$
for all $i = 1, \dots, n$.
\end{lemma}

%Lemma \ref{lemma:configuration}
\begin{proof}
We use the mirror of $\scC$
introduced by Khovanov, Seidel, and Thomas
\cite{Khovanov-Seidel, Seidel-Thomas}.
It is defined
%as the symplectic manifold obtained
as the deformation
$$
  V = \{ (x, y, z) \in \bC^3 \mid
   x y = z (z - 1) \cdots (z - n)
%   x y = z^{n+1} + a_{n-1} z^{n-1} + \dots + a_0
      \}
$$
of the $A_n$-singularity,
equipped with the symplectic form
given by restricting
the standard Euclidean K\"{a}hler form
on $\bC^3$.
The grading of $V$
is given by the holomorphic volume form
$
 \Omega = d y / y \wedge d z.
$
%The values of $(a_0, \dots, a_{n-2})$
%is irrelevant
%as long as they are outside the discriminant set
%so that
%$V$ is smooth,
%since the resulting $V$ is independent of it
%as a symplectic manifold.
There are $n$ graded Lagrangian submanifolds
$\{ \Ltilde_i \}_{i=1}^n$ of $V$,
where $L_i$ is the Lagrangian submanifold
coming from the line segment $c_i$
from $i-1$ to $i$ on the $z$-plane
using the conic fibration
$$
\begin{array}{rccc}
 \pi : & V & \to & \bC \\
 & \rotatebox{90}{$\in$} & & \rotatebox{90}{$\in$} \\
 & (x, y, z) & \mapsto & z
\end{array}
$$
in just the same way as
for $W$.
Let $\Fuk V$ be the Fukaya category of $V$
whose set of objects is
$\{ \Ltilde_i \}_{i=1}^n$.
Since
$(\Ltilde_i)_{i=1}^n$
is an $A_n$-configuration
generating $D^b \Fuk V$,
$D^b \Fuk V$ is equivalent to $\scC$
by the intrinsic formality
of Seidel and Thomas
\cite[lemma 4.21]{Seidel-Thomas}.
We fix this equivalence
by identifying $\scO_{C_i}(-1)$
with $\Ltilde_i$ for $i = 1, \dots, n$.

Now we prove Lemma \ref{lemma:configuration}
in characteristic two.
Then the lifting argument
in \S \ref{section:lifting}
shows that it holds
in any characteristic.
The assumption
$E_i \cong \Phi(\Ltilde_j)$
shows that
$E_i \in D^b \Fuk V$
is isomorphic
not only to a twisted complex of
graded Lagrangian submanifolds
but to an honest graded Lagrangian submanifold
constructed from a graded curve
$\widetilde{d}_i$ in $\bC$
whose endpoints lie in
$
 \Delta = \{ 0, \dots, n \}.
$
%Let $\Diff_0(\bC; \Delta)$
%be the group of diffeomorphisms of $\bC$
%which fixes $\Delta$ as a set.
%Then
%$\BG$ is isomorphic to
%$\pi_0(\Diff_0(\bC; \Delta))$,
%and
%we have to show that
%there is an element
%$g \in \Diff_0(\bC; \Delta)$
%such that $g(d_i) = c_i$.
Since
$(E_1, \dots, E_n)$
is an $A_n$-configuration,
the geometric intersection numbers
between the curves $d_i$ satisfy
$$
  I(d_i, d_j) =
   \begin{cases}
    \frac{1}{2} & \text{if } |i-j| = 1, \\
    0 & \text{if} \ |i-j| \geq 2,
   \end{cases}
$$
by Theorem \ref{theorem:Khovanov-Seidel}.
This shows that
$d_i$ intersects $d_{i-1}$ and $d_{i+1}$
only at its endpoints,
and does not intersect any other curve
$d_j$ with $|i-j|>1$.
Such a configuration of curves
can be taken
by the action of some
$g \in \Diff_0(\bC; \Delta)$
to the standard configuration
$(c_1, \dots, c_n)$
so that
$g(d_i) \simeq c_i$
for
$i = 1, \dots, n$.
Since
$\pi_0( \Diff_0(\bC; \Delta) ) \cong \BG$,
this suffices to show that
there is an element
$\Psi \in \Br(\scC)$
such that
$\Psi(E_i) \cong \Ltilde_i$.
\end{proof}

%%%%%%%%%%%%%%%%%%%%%%%%%%%%%%%%%%%%%%%%%%%%%%%%%%%%%%%%%%%
\subsection{The proof of the connectedness of $\Stab \scC$}

Here
we finish the proof of the connectedness of $\Stab \scC$.
Let $S_i$ be the simple $\palg$-module
corresponding to the idempotent $(i) \in \palg$
for $i = 1, \dots, n$,
and $\scA$ be the full extension-closed subcategory
of $\module_0 \palg$
containing $S_1, \dots, S_n$.
Then the standard $t$-structure of $D^b \module \palg$
induces a bounded $t$-structure of $\scC$ with heart $\scA$.
For any connected component $\Sigma$ of $\Stab \scC$,
take a stability condition
$\sigma \in \Sigma$
such that $Z(E) \not \in \bR $ for any spherical object $E \in \scC$.
Then $\sigma$ gives an $A_n$-configuration
$(E_1, \dots, E_n)$
by Lemma \ref{lemma:non-degenerate}.
Applying Lemma \ref{lemma:length} and Lemma \ref{lemma:configuration},
we can find an autoequivalence $\Psi \in \Br(\scC)$
such that
$\Psi(E_i) \cong S_i$.
Then the heart of $\Psi \sigma$ coincides with $\scA$
and hence $\sigma$
belongs to
the distinguished connected component
studied by Bridgeland in \cite{Bridgeland_SCKS}.

%%%%%%%%%%%%%%%%%%%%%%%%%%%%%%%%%%%%%%%%%%%%%%%%%%%%%%%%%%%%%%%%%%%%%%%%
%%%%%%%%%%%%%%%%%%%%%%%%%%%%%%%%%%%%%%%%%%%%%%%%%%%%%%%%%%%%%%%%%%%%%%%%
\appendix

\section{Appendix}

Let $\scD$ be the drived category
$D^b \coh_Z X$ of coherent sheaves
on the minimal resolution $X$ of an $A_n$-singularity
supported at the exceptional set $Z$
as in the body of this paper,
and $\Auteq^{FM}\scD$ be the subgroup of $\Auteq \scD$ 
consisting of integral functors.
%(See e.g. \cite[Introduction]{Ishii-Uehara_ADC}.)
A set of generators of $\Auteq^{FM}\scD$
is found in \cite[theorem 1.3]{Ishii-Uehara_ADC}.
Here we prove that every autoequivalence of $\scD$ is given by an integral functor,
\begin{equation}\label{equation:equiv}
 \Auteq^{FM}\scD= \Auteq \scD.
\end{equation} 
%We work over the complex number field $k=\bC$.
%This assumption is needed
%in the proof of Lemma \ref{lemma:FM}.

%%%%%%%%%%%%%%%%%%%%%%%%%%%%%%%%%%%%%%%%%%%%%%%%%%%%%%%%%%%
\subsection{Weak ample sequence and Fourier-Mukai transform}
Let $\scA$ be an abelian category
over a field $k$.
First we recall the definition of ample sequences in $\scA$,
introduced by Bondal and Orlov
\cite[\S 2]{Orlov_EDCKS},
\cite[appendix]{Bondal-Orlov}:

%%%
%%% 

\begin{definition}\label{definition:ample}
A collection $\{P_i\}_{i\in \bZ}$ of objects in $\scA$
is an {\em ample sequence}
if for any object $F\in \scA$,
there is an integer $N$ such that for any integer $i$
smaller than $N$,
the following conditions are satisfied:
\begin{enumerate}
\item
the canonical morphism $\Hom_{\scA}(P_i,F)\otimes P_i \to F$
is surjective,
\item
$\Ext_{\scA}^j(P_i,F)=0$ for any $j\ne 0$,
\item
$\Hom_{\scA}(F,P_i)=0$.
\end{enumerate}  
\end{definition}

When $X$ is a projective variety,
an ample line bundle $L$ on $X$ provides an example
$\{L^{\otimes i}\}_{i\in \bZ}$
of an ample sequence in $\coh X$.

Let $\{P_i\}_{i\in \bZ}$ be an ample sequence in $\scA$,
$D^b(\scA)$ be the bounded derived category of $\scA$,
and $\scB$ be the full subcategory of $D^b(\scA)$
consisting of $\{ P_i \}_{i \in \bZ}$.
The inclusion functor will be
denoted by $j : \scB \hookrightarrow D^b(\scA)$.
The following result is due to Bondal and Orlov:

%%%
%%%

\begin{proposition}%[\cite{Bondal-Orlov, Orlov_EDCKS}]
\label{proposition:greatOrlov}
Let $\Psi$ be an autoequivalence of $D^b(\scA)$
and assume that there is an isomorphism of functors
$g : j \stackrel{\sim}{\to} \Psi|_\scB$.
Then $g$ can be extended to an isomorphism
$id \stackrel{\sim}{\to} \Psi$
on the whole $D^b(\scA)$. 
\end{proposition}

The following lemma is proved in
\cite[claim 3.8.]{Ishii-Uehara_ADC}:

%%%
%%%
\begin{lemma}\label{lemma:FM}
Let $j':\module_0 \palg \hookrightarrow \scD$
be the inclusion
given by the McKay correspondence in \S \ref{section:McKay}.
Then for any autoequivalence
$\Phi \in \Auteq \scD$,
there is an autoequivalence
$\Psi \in \Auteq^{FM} \scD$
such that we have an isomorphism of functors
$
 h : j' \stackrel{\sim}{\to} \Psi \circ \Phi |_{\module_0 \palg}.
$
\end{lemma}

If $\module_0 \palg$ has an ample sequence,
then we can apply Proposition \ref{proposition:greatOrlov}
to obtain the equality (\ref{equation:equiv}).
But unfortunately, not all simple objects $F$ of $\module_0 \palg$ can satisfy the condition (iii).
We rather consider $\coh_Z X$.
Although
no sequence in $\coh_Z X$ satisfies
the conditions (i) and (ii) simultaneously for $F=\scO_x$,
we can weaken the condition (ii)
in the following way:

\begin{definition} \label{definition:weak_ample}
A sequence $\{ P_i \}_{i \in \bZ}$ of objects in $\scA$
is a {\em weak ample sequence}
if for any objects $F, G\in \scA$,
there is an integer $N$ such that for any $i<N$,
the condition (i) and (iii) in Definition \ref{definition:ample}
and the following condition (ii)$\,\!^\prime$ are satisfied;
\vspace{0mm}\\

(ii)$\,\!^\prime$
\begin{minipage}[t]{.9 \textwidth}
there are a natural number $l$ and a surjection 
$
 \varphi : P_i^{\oplus l} \to G
$
such that the pull-back
$$
 \varphi^* : \Ext_{\scA}^j(G,F) \to \Ext_{\scA}^j(P_i^{\oplus l},F)
$$
is the zero map for any $j\ne 0$.
\end{minipage}
\end{definition}
The proof of Proposition \ref{proposition:greatOrlov}
goes through also for weak ample sequences:

\begin{proposition} \label{proposition:weak}
Let $j : \scB \hookrightarrow D^b(\scA)$ be inclusion
of the full subcategory
consisting of a weak ample sequence of $\scA$.
Let $\Psi$ be an autoequivalence of $D^b(\scA)$
and assume that there is an isomorphism of functors
$g : j \stackrel{\sim}{\to} \Psi|_\scB$.
Then $g$ can be extended to an isomorphism
$id \stackrel{\sim}{\to} \Psi$
on the whole $D^b(\scA)$. 
\end{proposition}

%%%%%%%%%%%%%%%%%%%%%%%%%%%%%%%%%%%%%%%%

\subsection{A weak ample sequence in $\coh_Z X$}
Here we prove the equality (\ref{equation:equiv}).
First we show that $\coh_Z X$ has a weak ample sequence.
Consider $Z$ as the fundamental cycle
with its scheme structure.

\begin{lemma} \label{lemma:wacohz}
Let $H$ be an ample divisor on $X$
such that the divisor $H+Z$ is ample
and put
$$
 P_i:=\scO_{|i|Z}(iH)
$$
for $i \in \bZ$.
Then
$\{P_i\}_{i\in \bZ}$ 
is a weak ample sequence in $\coh_Z X$.
\end{lemma}

\begin{proof}
Let $F$ and $G$ be objects of $\coh_Z X$.
We show that there is an integer $N$
such that the conditions
(i), (ii)$\,\!^\prime$, and (iii) hold
for any integer $i$ smaller than $N$.

(i)
Let $N_1$ be a negative integer
such that $F$ is an $\scO_{- N_1 Z}$-module,
and $N_2$ be another integer such that
the canonical morphism
$$
 \Hom_{- N_1 Z}(\scO_{- N_1 Z}(N_2 H), F)
  \otimes \scO_{- N_1 Z}( N_2 H) \to F
$$
is surjective.
Then we can take $N = \min\{ N_1, N_2 \}$.

(iii)
Let $N$ be a negative integer such that
\begin{equation} \label{eq:homf}
 \Hom_X(F,\scO_Z(i(H + Z)) = 0 
\end{equation}
for any $i < N$.
The short exact sequences
$$
 0 \to \scO_{(-i-j-1)Z}(iH-(j+1)Z) 
   \to \scO_{(-i-j)Z}(iH-jZ) \to  \scO_Z(iH-jZ) \to 0
$$
for $j = 0, \ldots, - i - 2$
show that
$P_i$ has a filtration
whose subquotients consist of
$\scO_Z(i H - j Z)$ for $j = 0, \ldots, - i - 1$.
Since we have an inclusion
$$
 \scO_Z(i H - j Z) \hookrightarrow \scO_Z(i (H + Z))
$$
for $j \leq - i$,
(\ref{eq:homf}) implies
$$
 \Hom_X(F, \scO_{Z}(i H - j Z)) = 0
$$
and hence
$$
 \Hom_X(F, P_i) = 0.
$$

(ii)$\,\!^\prime$
Since $\Ext^2_X(P_i, F)=0$ for sufficiently small $i$ by (iii)
and the Serre duality,
we only need to check the case $j=1$.
Take a negative integer $N_1$
such that the canonical morphism
$$
 \varphi : \Hom_{X}(P_{i}, G) \otimes P_{i} \to G
$$
is surjective for any $i < N_1$.
We will prove that
there is an integer $N <  N_1$ such that
for any $i < N$ and any extension
\begin{equation} \label{eq:exactone}
 0 \to F \to E \to G \to 0
\end{equation}
of $G$ by $F$
corresponding to $e \in \Ext_X^1(G, F)$,
the extension
\begin{equation} \label{eq:exacttwo}
 0 \to F \to E' \to \Hom_{X}(P_{i}, G) \otimes P_{i} \to 0
\end{equation}
corresponding to
$$
 \varphi^*(e) = e \circ \varphi
  \in \Ext_X^1(\Hom_{X}(P_{i}, G) \otimes P_{i}, F)
$$
splits.

First note that if one takes a sufficiently small $i$
so that $F$ and $G$ are $\scO_{-i Z}$-modules,
then $E$ is an $\scO_{-2 i Z}$-module
for any $e \in \Ext_X^1(G, F)$.
Hence there is an integer $N_2 < N_1$ such that
$E$, $F$, and $G$ are $\scO_{- i Z}$-modules
for any $i < N_2$.
Since $H$ is ample,
there is an integer $N < N_2$ such that
$\Ext_{- i Z}^1(P_{i}, F) = 0$
for any $i < N$.
Take any $i < N$.
Then (\ref{eq:exacttwo}) is an exact sequence 
not only of $\scO_X$-modules
but also of $\scO_{- i Z}$-modules.
Since $\Ext_{- i Z}^1(P_{i}, F) = 0$,
it splits as an exact sequence
of $\scO_{- i Z}$-modules,
and hence as an exact sequence of $\scO_X$-modules.
\end{proof}

The above weak ample sequence of $\coh_Z X$
belongs also to $\module_0 \palg [-1]$:

\begin{lemma} \label{lemma:per}
If $i < 0$, then $P_i$ in Lemma \ref{lemma:wacohz}
belongs to $\Per{-1}(X/Y)[-1]$,
which is identified with $\module \palg[-1]$
by the McKay correspondence in \S \ref{section:McKay}.
\end{lemma}

\begin{proof}
It follows from the definition of $\Per{-1}(X/Y)$
that a coherent sheaf $F$ on $X$
belongs to $\Per{-1}(X/Y)[-1]$
if and only if $f_* F = 0$.
As in the proof of (iii) in Lemma \ref{lemma:wacohz},
we have $f_* P_i = 0$ for $i < 0$
and hence $P_i \in \Per{-1}(X/Y)[-1]$.
\end{proof}

Now we can prove the equality (\ref{equation:equiv}).
Let $\Phi$ be an autoequivalence of $\scD$.
Then there is another autoequivalence
$\Psi \in \Auteq^{FM} \scD$
such that we have an isomorphism of functors
$$
 h : j' \stackrel{\sim}{\to} \Psi \circ \Phi |_{\module_0 \palg}
$$
by Lemma \ref{lemma:FM}.
Let $j : \scB \hookrightarrow \scD$ be
the inclusion of the full subcategory
consisting of the weak ample sequence $\{ P_i \}_{i < 0}$
given in Lemma \ref{lemma:wacohz}.
Since
$
 \scB \subset \module_0 \palg [-1]
$
by Lemma \ref{lemma:per},
$h$ induces an isomorphism
$$
 g : j \stackrel{\sim}{\to} \Psi \circ \Phi |_\scB
$$
of functors from $\scB$ to $\scD$.
Since $\{ P_i \}_{i \in \bZ}$ is a weak ample sequence
in $\coh_Z X$ by Lemma \ref{lemma:wacohz} and
$\scD = D^b \coh_Z X$,
$g$ can be extended to an isomorphism of functors
on the whole $\scD$ by Proposition \ref{proposition:weak}.

%\bibliographystyle{plain}
%\bibliography{bibs}

\newcommand{\noop}[1]{}\def\cprime{$'$} \def\cprime{$'$}
  \def\cftil#1{\ifmmode\setbox7\hbox{$\accent"5E#1$}\else
  \setbox7\hbox{\accent"5E#1}\penalty 10000\relax\fi\raise 1\ht7
  \hbox{\lower1.15ex\hbox to 1\wd7{\hss\accent"7E\hss}}\penalty 10000
  \hskip-1\wd7\penalty 10000\box7} \def\cydot{\leavevmode\raise.4ex\hbox{.}}
  \def\cprime{$'$} \def\cprime{$'$} \def\cprime{$'$}

\noindent
Akira Ishii

Department of Mathematics,
Graduate School of Science,
Hiroshima University,
1-3-1 Kagamiyama,
Higashi-Hiroshima,
739-8526,
Japan

{\em e-mail address}\ : \ akira@math.sci.hiroshima-u.ac.jp

\ \\

\noindent
Kazushi Ueda

Department of Mathematics,
Graduate School of Science,
Osaka University,
Machikaneyama 1-1,
Toyonaka,
Osaka,
560-0043,
Japan.

{\em e-mail address}\ : \  kazushi@math.sci.osaka-u.ac.jp

\ \\

\noindent
Hokuto Uehara

Department of Mathematics and Information Sciences,
Tokyo Metropolitan University,
1-1 Minami-Ohsawa,
Hachioji-shi,
Tokyo,
192-0397,
Japan 

{\em e-mail address}\ : \  hokuto@tmu.ac.jp

\ \\

\end{document}